\documentclass[10pt]{article}
\usepackage[utf8]{inputenc}
\usepackage{latexsym,amssymb,amsmath,url,authblk}

\def\N{\mathbb{N}}
\def\Q{\mathbb{Q}}
\def\Z{\mathbb{Z}}
\def\R{\mathbb{R}}

\def\proof{\par\noindent{\em Proof. }}
\def\eproof{\hfill{$\Box$}\bigskip}

\def\ds{\dots}
\def\sus{\subset}

\def\de{\delta}
\def\la{\lambda}
\def\cc{\colon}

\newtheorem{thm}{Theorem}[section]
\newtheorem{prop}[thm]{Proposition}
\newtheorem{cor}[thm]{Corollary}
\newtheorem{lem}[thm]{Lemma}

\newtheorem{defi}[thm]{Definition}

\title{Countable real analysis}
\author{Martin Klazar}
\date{\today}

\begin{document}

\maketitle
\centerline{{\em Dedicated to the memory of my father Ji\v r\'\i\ Klazar (1934--2022)}}

\begin{abstract}
HMC sets are hereditarily at most countable sets. We rework 
a~substantial part of 
univariate real analysis in a~form in which only HMC real functions are 
used. In such countable real analysis we carry out 
Hilbert's proof of transcendence of the number $\mathrm{e}$. We also construct a~uniformly continuous 
function $f:[0,1]\cap\Q\to\R$ such that $f'=1$ on $[0,1]\cap\Q$ and 
$\lim_{\substack{a\to1/\sqrt{2}\\a\in\Q}}f(a)=\frac{1}{\sqrt{2}}>f(b)$
for every $b\in[0,1]\cap\Q$. 
\end{abstract}

\tableofcontents

\section{Introduction}\label{sec_intro}

A~set is {\em at most countable} if it is finite or countable. In the former case 
it is in bijection with a~natural number 
$n=\{0,1,\ds,n-1\}$, $n\in\omega$, and in the latter case with the set $\omega=\{0,1,\ds\}$. A~set is {\em uncountable} if it is not at most countable. A~set 
$x$ is {\em hereditarily at most countable}, abbreviated HMC, if for every index $n\in\omega$ and every chain of sets
$$
x_n\in x_{n-1}\in\ds\in x_0=x\;,
$$
the set $x_n$ is at most countable. See \cite{jech82} for a~closely related notion. Due 
to the axiom of foundation every chain of sets 
$x_0\ni x_1\ni\ds$ is finite.  See \cite{jech,schi,soch} for set-theoretical terminology and notions.

Are uncountable sets/functions indispensable in (univariate real) analysis? Many people think so. In this article we argue to the contrary; in 
Section~\ref{sec_HilbProof} we carry out, by means of only HMC real
functions, Hilbert's proof \cite{hilb} of transcendence of the 
number~$\mathrm{e}$. Hilbert's proof is based on an integral identity due 
to L.~Euler: for every $n\in\omega$, 
$$
\int_0^{+\infty}x^n\mathrm{e}^{-x}\,\mathrm{dx}=n!\ \ (=1\cdot2\cdot\ldots\cdot n)\;.
$$
We replace the uncountable integrands with their restrictions to $\{a\in\Q:\;a\ge0\}$ and adapt improper integrals for such countable real functions. 
By a~real function we mean any function of the type
$f\cc M\to\R$ with $M\sus\R$. Assuming that real numbers are Dedekind cuts, such a~function
is HMC if and only if it is at most countable. Abusing terminology a~little, instead of HMC real 
functions we usually speak shortly of countable functions.

In Section~\ref{sec_inteRealF} we rework results in univariate real 
analysis needed for Hilbert's proof so that they use only countable 
functions, see the above table of contents. 
Subsections~\ref{subsec_sequences} and \ref{subsec_series} contain 
standard 
material on real sequences and real series which are trivially 
countable functions. Our treatment of real numbers in 
Subsection~\ref{subsec_reals} is of some interest. Proofs in these 
three subsections are included 
for completeness. In Subsection~\ref{subsec_HMCfun} we 
make a~sharp turn away from the main road of real analysis and start to develop 
a~theory of functions $f\cc M\to\R$ with $M\sus\Q$. Now the proofs are completely 
new. Results in real analysis are usually not too hard to adapt for countable 
functions, possibly with the exception of Theorem~\ref{thm_GlobEx}.

Our adaptation of real analysis rests on two innovations. {\bf We only use 
uniformly continuous (UC) functions.} Uniform continuity takes over the role of 
compactness of definition domains. UC functions suffice because real 
functions in analytic number theory are often locally UC. Euler's identity is derived by means of the 
fundamental theorem of analysis, which in turn follows from 
Lagrange's mean value theorem. This theorem follows from Rolle's 
theorem. Rolle's theorem follows from  
the {\em min-max principle}\,---\,any continuous 
function defined on a~compact set attains minimum and maximum\,---\,and 
from the {\em vanishing derivative principle} (VDP)\,---\,any function 
with a~nonzero derivative at a~two-sided limit point of the 
definition domain does not have  local extreme at the point. It 
is this latter VDP that is not straightforward to adapt.

The set $[0,1]_{\Q}=[0,1]\cap\Q$ is bounded but not compact, and it is no wonder that the countable function $f(a)=\frac{1}{a-1/\sqrt{2}}\cc [0,1]_{\Q}\to\R$ is continuous and unbounded. But it is not UC.
We recover the min-max principle for countable functions by means of  uniform continuity in  
Theorem~\ref{thm_attaGlobExtr}:
\begin{quote}
{\em If $M\sus\Q$ is bounded and $f\cc M\to\R$ is {\em UC}, then there exist real numbers 
$y,y'\in\overline{M}$ such that
$f(y)\le f(a)\le f(y')$ for every $a\in M$.} 
\end{quote}
Thus $f$ ``attains'' at the points $y$ and $y'$ in the closure of $M$ minimum and maximum, even if these points are not in $M$. The ``values'' $f(y)$ and $f(y')$ of $f$ are obtained by means of Theorem~\ref{thm_extending}: 
\begin{quote}
{\em If $M\sus\Q$ and $f\cc M\to\R$ is {\em UC} then $f$ has a~unique continuous 
extension $f\cc M\cup\{x\}\to\R$ to any 
point $x\in\overline{M}$ in the closure  of the definition domain $M$.} 
\end{quote}

{\bf We adapt VDP by means of the second innovation, the uniform derivative which is introduced in Definition~\ref{def_unifDeri}:}

\begin{quote}
{\em If $M\sus\Q$, $f\cc 
M\to\R$ and $D(f)\sus M$ is the
definition domain of $f'$, then we call $f'$ uniform if for every $k\in\N$ there is an $l\in\N$ such that}
$${\textstyle
a\in D(f)\wedge b\in M\wedge0<|b-a|
\le \frac{1}{l}\Rightarrow\big|\frac{f(b)-f(a)}{b-a}-f'(a)\big|\le \frac{1}{k}\,.
}
$$
\end{quote}
The adapted VDP emerges in Theorem~\ref{thm_GlobEx}:
\begin{quote}
{\em 
If $u<x<v$ are real numbers, $I=[u,v]_{\Q}$,  
$f,f'\cc I\to\R$, $f'$ is uniform and {\em UC} and $f'(x)\ne0$, then 
there exist $b,b'\in I$ such that
$$
f(b)<f(x)<f(b')\,.
$$
Hence $f(x)$ is not a~global extreme of $f$.
}    
\end{quote}
By Proposition~\ref{prop_UCandDer}, $f$ is UC as well and the values 
$f'(x)$ and $f(x)$ are understood in the sense of extension provided 
by Theorem~\ref{thm_extending}. With Theorems~\ref{thm_GlobEx} and 
\ref{thm_attaGlobExtr} it is easy to adapt Rolle's theorem and eventually get countable  
Euler's identity. On the way we adapt improper integration and the 
exponential function. Then it is easy to rewrite Hilbert's proof in terms of countable functions.

Recall that $[0,1]_{\Q}=[0,1]\cap\Q$. In 
Section~\ref{sec_stranFun} in Theorem~\ref{thm_nonAvoidUniDer} we demonstrate by an example that the assumption of uniform derivatives in 
Theorem~\ref{thm_GlobEx} is substantial. We construct 
a~function $f\cc[0,1]_{\Q}\to\R$ such that 
\begin{itemize}
    \item $f$ is UC (uniformly continuous).
    \item $f(\frac{1}{\sqrt{2}})\equiv\lim_{\substack{a\to\frac{1}{\sqrt{2}}\\ a\in\Q}}f(a)=\frac{1}{\sqrt{2}}$ and $\frac{1}{\sqrt{2}}>f(b)$ for every $b\in[0,1]_{\Q}$.
    \item $f'(b)=1$ for every $b\in[0,1]_{\Q}$.
\end{itemize}
The graph of $f$ consists of 
countably many straight segments with slope~$1$: the long one 
$\{(a,a)\in\Q^2\;|\;0\le a<\frac{1}{\sqrt{2}}\}$ and 
countably many short ones lying between the lines
$x=\frac{1}{\sqrt{2}}$ and $x=1$ and below the line
$y=\frac{1}{\sqrt{2}}$. The segments are mutually separated by 
vertical lines with irrational $x$-coordinates. Hence $f'=1$ on 
$[0,1]_{\Q}$. Since $f$ refutes the conclusion of 
Theorem~\ref{thm_GlobEx} (the adapted VDP), the constant 
derivative $f'=1$ is not uniform on $[0,1]_{\Q}$. 

The {\em countable univariate real analysis} which we develop in 
Section~\ref{sec_inteRealF} uses only the following non-HMC sets, which are in general uncountable.
\begin{enumerate}
\item The set $K$ of cuts on $\Q$ introduced in Definition~\ref{def_reals}. 
This is the set of real numbers. We also use real intervals $I$ and, more generally, 
closures $\overline{M}$ in $K$ of sets $M\sus\Q$.
\item The three sets $+,\cdot\cc K\times K\to K$ and $<\,\sus K\times K$, introduced in Definition~\ref{def_realArithm}. These 
are the addition, multiplication and ordering of real numbers.
\item The sets $\Q^{\N}$ and $C$ of rational sequences and of Cauchy rational sequences. 
\item The sets $\mathcal{F}(M)$,  $M\sus\Q$, of functions $f\cc M\to K$. 
Note that every function 
$f\in\mathcal{F}(M)$ is an HMC set.
\end{enumerate} 
Proposition~\ref{prop_NthPoRin} allows in $R^{\N}$ any set $R$, but it is applied 
only with $R=\Q$. We use the non-HMC sets in 1--4 only as notational conveniences; all 
can be easily formulated away by more complicated notation. For example, instead of 
$x\in K$ we can write 
$\exists\,y\,\big(x=y\wedge y\sus\Q\wedge K(y)\big)$, with $K(y)$ being the definition of 
a~cut, and similarly for $I$, $\overline{M}$, $+$, $\cdot$, $<$, 
$\Q^{\N}$, $C$ and $\mathcal{F}(M)$, $M\sus\Q$. This is possible due to the 
fact that the main objects of definitions and claims (theorems, 
propositions, $\ds$) of countable real analysis are not these sets but their 
elements, especially elements of the sets $\mathcal{F}(M)$, $M\sus\Q$, which are always HMC sets. Since we eliminate from 
(univariate real) analysis all uncountable functions by working only with functions of the
type $f:M\to\R$, $M\sus\Q$, does it mean that we eliminate from analysis all  
uncountable sets? In practice yes, in theory no. In practice we are playing 
in our playground only with functions $f\in\mathcal{F}(M)$, $M\sus\Q$, which are all HMC sets. In theory, beyond the fence 
of the playground in the forest, uncountable sets are still lurking. The set $\omega=\{0,1,
\ds\}$, on which everything is built, is defined in ZFC by means of the axiom of infinity as 
an intersection of (all) inductive sets. These are of completely undetermined nature and some 
of them are uncountable.  In every formal proof of every claim (theorem) involving 
$\omega$ this definition is always present. 

{\bf The main contributions of our article are 
three:} (i) the body of definitions and proven claims in the next Section~\ref{sec_inteRealF}, (ii) 
Theorem~\ref{thm_nonAvoidUniDer} (UC function $f\cc[0,1]_{\Q}\to\R$ with 
constant derivative $f'=1$ and strict
global maximum $f(\frac{1}{\sqrt{2}})
=\frac{1}{\sqrt{2}}$), and (iii) the proof of transcendence of the number $\mathrm{e}$ 
in Theorem~\ref{thm_Hermite} by means of HMC real functions only. 

In the final Section~\ref{sec_conclRem} we give several remarks, one concerning 
the definition of planar graphs in finite combinatorics. It relies on uncountable sets! 

\section{Countable real analysis}\label{sec_inteRealF}

We use $\equiv$ as a~defining equality, thus in $x\equiv y$ we define the new symbol $x$ by the 
already known expression $y$; 
exceptionally $x$ and $y$ may exchange their roles. 
We denote sets both by capital 
and small letters, capital letters do not mean classes. By $\sus$ we denote  non-strict inclusion of sets. Let $X$ and $Y$ be sets.
By writing $f\cc X\to Y$ we say that $f$ is a~{\em function (map)} from 
$X$ to $Y$ --- $f\sus X\times Y$ and for every $x\in X$ there is a~unique $y\in Y$ such that 
$(x,y)\in f$, symbolically $f(x)=y$. Here
$(x,y)=\{\{y,x\},\{x\}\}$ is Kuratowski's ordered pair.
We call $X$ the {\em definition domain} of $f$. Let $Z$ be any set.
The {\em restriction} of $f\cc X\to Y$ to $Z$ is the function 
$f\,|\,Z\cc Z\cap X\to Y$ with values 
$$
(f\,|\,Z)(x)\equiv f(x),\ \ x\in Z\cap X\;. 
$$
Instead of $f\,|\,Z$ we often write just $f$. The {\em image} of $Z$ by $f$ is the set
$$
f[Z]\equiv\{f(x):\;x\in Z\cap X\}\ \ (\sus Y)\,.
$$
The {\em preimage} of $Z$ by $f$ is the set
$$
f^{-1}[Z]\equiv\{x\in X:\;f(x)\in Z\}\ \ (\sus X)\,.
$$

We use the sets
$$
\omega\equiv\{0,\,1,\,\ds\}\,\text{ and }\,
\N\equiv\omega\setminus\{0\}\;.
$$
Here $0\equiv\emptyset$, $1\equiv\{0\}$, $2\equiv\{0,1\}$ and so on. The elements of $\omega$ are 
called {\em natural numbers}. On $\omega$ we have the standard arithmetic operations 
$+$ and $\cdot$. They are associative and commutative and $\cdot$ is distributive to $+$. The element 
$0$ (resp. $1$) is neutral to $+$ (resp. $\cdot$). 
So $(\omega,0,1,+.\cdot)$ is a~semiring. The (binary) relation
$(\omega,<)\equiv(\omega,\in)$ is a~well ordering --- it is a~linear order on $\omega$ with the property that every nonempty subset of $\omega$ has the least element. If $A$ is a~set, 
a~{\em sequence} in $A$ is a~function $x\cc\N\to A$. We set $x_m\equiv x(m)$ and write $(x_n)\sus A$ to invoke a~sequence in $A$.
 
$$
\Z\equiv
\{\ds,\,-2,\,-1,\,0,\,1,\,2,\,\ds\}=(-\N)\cup\omega
$$ 
is the set of {\em integers}. Here
$-\N\equiv\{-n:\;n\in\N\}$ and $-n\equiv(n,0)$.
With respect to the standard addition $+$, multiplication $\cdot$ and linear ordering $<$, the integers form an ordered integral domain. For $n\in\N$ we define 
$$
[n]\equiv\{1,\,2,\,\ds,\,n\}\,,
$$
and set $[0]\equiv\emptyset$. The symbols $k$, $l$, $m$, $n$, $n_0$, $n_1$, $\ds$, $n_1'$, 
$n_2'$, $\ds$ denote elements of $\N$, and $k$, $l$, $m$ and $n$ range also in $\Z$. 

The elements of the set 
$$
\Q\equiv
\{[m/n]_{\sim}:\;m,\,n\in\Z,\,n\ne 0\}\ \ \text{(here $m/n\equiv(m,\,n)$)}
$$ 
are blocks of the equivalence relation $\sim$ on $\Z\times(\Z\setminus\{0\})$, where
$$
m/n\sim k/l\equiv m\cdot l=n\cdot k\;.
$$
These equivalence blocks are called {\em fractions} or {\em rational numbers}. Instead of $[m/n]_{\sim}$ we write usually just  $m/n$
or $\frac{m}{n}$. The algebraic structure
$${\textstyle
\Q_{\mathrm{OF}}\equiv\big(\Q,\,\frac{0}{1},\,\frac{1}{1},\,+,\,\cdot,\,<\big)
}
$$ 
with the standard arithmetic operations $+$ and $\cdot$ and the standard linear order $<$ is an ordered field. We denote fractions by the letters $a$, $b$, 
$c$ and $d$, possibly enriched with indices 
and/or primes. We introduce {\em rational intervals}.

\begin{defi}\label{def_ratioInter}
$J\sus\Q$ is a~rational interval if for any $a<b<c$ with 
$a,c\in J$ also $b\in J$. $J$ is nontrivial 
if $J\ne\emptyset,\{a\}$. If $I$ 
is a~real interval, then $I_{\Q}\equiv I\cap\Q$ and is 
a~rational interval. More precisely, $I_{\Q}\equiv I\cap E[\Q]$, see below.
\end{defi}

The map
$$
\Z\ni n\mapsto [n/1]_{\sim}\in\Q
$$
embeds integers in the field $\Q$ as a~subring. Instead of $n/1$ we 
usually write just $n$. We  simplify other rational expressions 
in a~similar way. 
For example, if $a\in\Q$ then $\frac{a}{2}$ is really $\frac{a}{[2/1]_{\sim}}=a\cdot[1/2]_{\sim}$.

\subsection{Real numbers as cuts}\label{subsec_reals}

We start from the real numbers. We define them as Dedekind cuts \cite{dede}, then they are HMC sets. We do not want to use their 
definition as Cantor's equivalence blocks because each of these blocks is uncountable.

\begin{defi}\label{def_reals}
A~cut is a~set $X\sus\Q$ such that {\em (i)} $X,\Q\setminus X\ne\emptyset$, {\em (ii)} for any $a\in\Q$, $b\in X$ with $a<b$ also $a\in X$
and {\em (iii)} $X$ has no largest element. We denote the set of cuts by $K$. We call the elements of $K$ also real numbers.
\end{defi}
Let $X$ be a~cut. If $a\in X$, $b\in\Q\setminus X$ then
$a<b$, and if $a\in\Q\setminus X$, $b\in\Q$ $\wedge$ $a<b$ then $b\in\Q\setminus X$. The injective map
$$
E\cc\Q\to K,\ E(m/n)\equiv
\{a\in\Q:\;a<m/n\}\,,
$$
embeds the field $\Q_{\mathrm{OF}}$ as an ordered subfield of the ordered field 
$$
\R\equiv(K,\,0_K,\,1_K,\,+,\,\cdot,\,<)\,.
$$
of real numbers. That $\R$ is an ordered field will be proven in 
Theorem~\ref{thm_KisOF}. So any 
rational number
$\frac{m}{n}$ is also a~real number $E(m/n)$. We usually do 
not distinguish between $a\in\Q$ and $E(a)\in K$. For 
example, for $x\in K$ and $k\in\N$ we write just $x\le\frac{1}{k}$ 
instead of $x\le E(\frac{1}{k})$ or $x\le 1_K/E(k)$. We introduce {\em 
irrational (real) numbers}.

\begin{defi}\label{def_irrNum}
$X\in K$ is an irrational number if $X\not\in E[\Q]$.   
\end{defi}
Starting from the next subsection we denote real numbers by the letters $u$, $v$, $w$, $x$, 
$y$ and $z$, possibly enriched with indices and/or primes. In this subsection we use capital letters $W$, $X$,  
$Y$ and $Z$.

In \cite{klaz_dedeCut} we devised {\em generating sequences} and with their help we endow $K$ with 
arithmetic operations. We say that $(a_n)\sus\Q$ is a~{\em Cauchy sequence} if for every $k$ there is an $n_0$ such that
$${\textstyle
m,\,n\ge n_0\Rightarrow|a_m-a_n|\le\frac{1}{k}\,.
}
$$

\begin{lem}\label{lem_onCauSeq}
If $(a_n)\sus\Q$ and $b$ are fractions such that $a_1\ge a_2\ge\ds\ge a_n\ge\ds\ge b$, then $(a_n)$ is a~Cauchy sequence. The same holds if all inequalities $\ge$ are reversed.
\end{lem}
\proof
We prove just the former claim, proof for 
the latter follows by reversing inequalities. It follows from the assumptions that for every 
$m<n$ we have $0\le a_m-a_n\le a_1-b$. Suppose that $(a_n)$ is not Cauchy. Then there is a~$k$ and 
natural numbers $n_1<n_2<\ds$ such that $a_{n_m}-a_{n_{m+1}}\ge\frac{1}{k}$ for every 
$m$. For any large $m$ such that $\frac{m}{k}>a_1-b$ we get the contradiction that
$${\textstyle
a_{n_1}-a_{n_{m+1}}=\sum_{j=1}^m(a_{n_j}-
a_{n_{j+1}})\ge\frac{m}{k}>a_1-b\,.}
$$
\eproof

\noindent
This lemma is an ur-theorem on limits of monotone sequences.   

\begin{defi}\label{def_C}
$C$ denotes the set of Cauchy sequences.  
\end{defi}
We define an equivalence relation $\sim$ on $\Q^{\N}$ by setting $(a_n)\sim(b_n)$ iff
for every $k$ there is an $n_0$ such that
$n\ge n_0\Rightarrow|a_n-b_n|\le\frac{1}{k}$.
If $(a_n)\sim(0,0,\ds)$, we write that $a_n\to0$. We relate Cauchy 
sequences and cuts.

\begin{defi}\label{def_G}
For $(a_n)\sus\Q$ let
$$
G(a_n)\equiv
\{a\in\Q:\;\exists\,n_0\,\big(n\ge n_0\Rightarrow a\le a_n\big)\}'\ \ (\sus\Q)\,.
$$
The prime indicates that the largest element, if present, is deleted from the set. 
We say that the sequence $(a_n)$ generates the set $G(a_n)$.
\end{defi}
Recall that $C$ is the set of Cauchy sequences and $K$ is the set of cuts.

\begin{thm}\label{thm_CandK}
The following hold.
\begin{enumerate}
    \item For every $(a_n)\in C$ the set $G(a_n)$ is a~cut.
    \item For every cut $X$ there is an $(a_n)\in C$ such that $X=G(a_n)$.
    \item If $(a_n),(b_n)\in C$ then $(a_n)\sim(b_n)$ if and only if $G(a_n)=G(b_n)$.
\end{enumerate}
Hence the map $G\cc C/\!\sim\,\to K$, $[(a_n)]_{\sim}\mapsto G(a_n)$, is a~bijection. 
\end{thm}
\proof
1. Let $(a_n)\in C$ (it suffices if $(a_n)$ is bounded) and $X\equiv G(a_n)$. The set $\{a_n:\;n\in\N\}$ ($\sus\Q$) is 
bounded, hence there are fractions lying below and above it. So $X,\Q\setminus 
X\ne\emptyset$ and condition (i) in Definition~\ref{def_reals} holds. If 
$b\in X$ then $b\le a_n$ for any large~$n$. It holds also for any 
fraction $a<b$. Thus $a\in X$ and condition (ii) holds. We show that (iii) holds. Suppose for the 
contrary that $a=\max(X)$. Let $S\equiv
\{\cdots\}$ be the above displayed set in the definition of  
$G(a_n)$. If $\max(S)$ does not exist then $X=S$. This is impossible because $X$ does have a~maximum. If $s\equiv\max(S)$ exists then $a<s$. Then every fraction between $a$ and $s$ belongs to
$X$, which contradicts the maximality of $a$. Both cases lead to a~contradiction. So (iii) holds and $X$ is a~cut.

2. Let $X\in K$. We define sequences $(a_n)\sus X$ and 
$(b_n)\sus\Q\setminus X$ (hence always $a_m<b_n$) such that $b_n-
a_n\to0$. We begin with arbitrary $a_1\in X$ and $b_1\in\Q\setminus 
X$. If $a_1,\ds,a_n,b_1,\ds,b_n$ are already defined, we set 
$c\equiv\frac{a_n+b_n}{2}$. If $c\in X$ then
$a_{n+1}\equiv c$ and $b_{n+1}\equiv b_n$. If $c\not\in X$ then $a_{n+1}\equiv a_n$ and $b_{n+1}\equiv c$.
By Lemma~\ref{lem_onCauSeq}, $(a_n),(b_n)\in C$. We show that $X=G(a_n)$. Since $(a_n)\sus X$ and $X$ is 
a~cut, it is clear that $X\supset G(a_n)$. Let $a\in X$ be
arbitrary. Then $a<a'<b$ for some $a'\in X$ and every $b\in\Q\setminus X$. Hence, since $(b_n)\sus 
\Q\setminus X$ and $b_n-a_n\to 0$, for every large $n$ one has that $a<a_n$ and therefore $a\in G(a_n)$. 
Thus $X\sus G(a_n)$ and $X=G(a_n)$.

3. Let $(a_n),(b_n)\in C$. Suppose that $(a_n)\sim(b_n)$ and $a\in G(a_n)$. 
Since $G(a_n)$ is a~cut, $a<a'<a_n$ for some $a'\in G(a_n)$ and every large $n$. Since $(a_n)\sim(b_n)$,
for every large $n$ one has that $a<b_n$ and $a\in G(b_n)$. Thus $G(a_n)\sus G(b_n)$. 
The same argument gives the opposite inclusion and $G(a_n)=G(b_n)$.

Suppose that $G(a_n)=G(b_n)\equiv X$. Let a~$k$ be given. The two sequences in part~2 show that
we can take fractions $a<b$ such that $a\in X$, $b\in\Q\setminus X$ and $b-a\le\frac{1}{2k}$. Since $(b_n)\in C$, for every large $m$ and $n$
it holds that $|b_m-b_n|\le\frac{1}{2k}$. Since $a\in G(a_n)$ and $b\not\in G(b_n)$, for every large $n$ 
it holds that $a_n\ge a$ and for infinitely many $n$ we have that 
$b_n<b$. Thus for every large $n$,
$${\textstyle
b_n-a_n\le b-a+\frac{1}{2k}\le\frac{1}{k}\,.
}
$$
Similar argument, with $(a_n)$ and $(b_n)$ exchanged, shows the 
opposite inequality that for every large $n$ one has that
$a_n-b_n\le\frac{1}{k}$. Thus for every large $n$ we have that $|a_n-
b_n|\le\frac{1}{k}$ and $(a_n)\sim(b_n)$.
\eproof

\noindent
We did not find this theorem in the literature \cite{buko,stil} 
on real numbers. The map
$G\cc C/\!\!\sim\,\to K$ transfers the straightforward arithmetic of 
Cantor's blocks to Dedekind cuts. The arithmetic of cuts is usually 
introduced directly, but then it is a~bit cumbersome. 

We define zero, identity, addition, multiplication, and ordering in $K$.

\begin{defi}\label{def_realArithm}
$0_K\equiv G(\frac{0}{1},\frac{0}{1},\ds)=E(\frac{0}{1})$ and $1_K\equiv G(\frac{1}{1},\frac{1}{1},\,\ds)=E(\frac{1}{1})$. 
For cuts $X,Y\in K$ with $X=G(a_n)$ and $Y=G(b_n)$ (we use part~2 of Theorem~\ref{thm_CandK}) we define their sum and product as 
$$
X+Y\equiv G(a_n+b_n)\,\text{ and }\,X\cdot Y\equiv G(a_nb_n)\,.
$$
We define 
$${\textstyle
X<Y\iff\exists\,k,\,n_0\,\big(
n\ge n_0\Rightarrow a_n\le b_n-\frac{1}{k}\big)\,.
}
$$ 
\end{defi}
We show that these definitions are correct. 

\begin{prop}\label{prop_correctArtithm}
{\em 1. }$0_K,1_K\in K$. {\em 2. }If $(a_n),(b_n)\in C$ then $(a_n+b_n)\in C$ and $(a_n b_n)\in C$. {\em 3. }The sum and product of $X=G(a_n)$ and $Y=G(b_n)$ and their comparison by $<$ do not depend on the generating sequences $(a_n)$ and $(b_n)$. 
\end{prop}
\proof
1. Since $(\frac{0}{1},\frac{0}{1},\ds)$ and $(\frac{1}{1},\frac{1}{1},\ds)$ are 
are bounded, $0_K$ and $1_K$ are cuts by 1~of Theorem~\ref{thm_CandK}.

2. Let $(a_n),(b_n)\in C$. It is easy to see that $(a_n+b_n)\in C$. We treat $(a_nb_n)$ in detail. 
Since $(a_n)$ and $(b_n)$ are Cauchy, they are bounded and there is an $l$ such
that for every $n$ one has that $|a_n|,|b_n|\le l$. Now for a~given $k$ for every large $m$ and $n$ it holds that  
$|a_m-a_n|,|b_m-b_n|\le\frac{1}{k}$. Then for the same $m$ and $n$ we have that
$${\textstyle
|a_mb_m-a_nb_n|\le|a_m|\cdot|b_m-b_n|+|a_m-a_n|\cdot|b_n|\le\frac{2l}{k}\,.
}
$$
Hence $(a_nb_n)\in C$.  

3. Suppose that $(a_n)$, $(b_n)$, $(a_n')$ and $(b_n')$ are in $C$ 
and that $X=G(a_n)=G(a_n')$ and $Y=G(b_n)=G(b_n')$. 
By part~3 of Theorem~\ref{thm_CandK}, $(a_n)\sim(a_n')$
and $(b_n)\sim(b_n')$. As in part~2 we show that 
$(a_n+b_n)\sim(a_n'+b_n')$ and $(a_nb_n)\sim(a_n'b_n')$. Thus 
again by part~3 of Theorem~\ref{thm_CandK},
$$
G(a_n+b_n)=G(a_n'+b_n')\,\text{ and }\,G(a_nb_n)=G(a_n'b_n')\;.
$$
So the cuts $X+Y$ and $X\cdot Y$ do not depend on the choice of generating sequences. We show it  also for comparison. Let $X<Y$. It means that for some $k$ and every large 
$n$ one has that $a_n\le b_n-\frac{1}{k}$ . Since $(a_n)\sim(a_n')$
and $(b_n)\sim(b_n')$, for every large $n$ one has that 
$|a_n-a_n'|,|b_n-b_n'|\le\frac{1}{4k}$. Then for every large $n$ we have that
$${\textstyle
a_n'\le a_n+\frac{1}{4k}\le b_n-\frac{1}{k}+\frac{1}{4k}\le b_n'-\frac{1}{k}+\frac{1}{4k}+\frac{1}{4k}=b_n'-\frac{1}{2k}\,.
}
$$
Thus $X<Y$ holds irrespective of the choice of generating sequences.
\eproof

We express the relation $<$ on $K$ in terms of $\sus$. The notation $X\le Y$ means that $X<Y$ or $X=Y$. 

\begin{prop}\label{prop_lesVerInc}
Let $X,Y\in K$. Then  
$X\le Y$ if and only if $X\sus Y$. 
\end{prop}
\proof
Let $X=G(a_n)$ and $Y=G(b_n)$ with $(a_n),(b_n)\in C$. Suppose that 
$X\le Y$. If $X=Y$ then also $X\sus Y$. Let $X<Y$. Then for some $k$ and 
every large $n$ one has that $a_n\le b_n-\frac{1}{k}$. If $a\in X$ then 
$a\le a_n$ for every large $n$, hence also $a\le b_n$ for every 
large $n$ and $a\in G(b_n)=Y$. So $X\sus Y$. 

Suppose that $X\sus Y$. If $X=Y$, there is nothing to prove. We therefore assume that $X\ne Y$
and take fractions $b'<b$ such that $b',b\in Y\setminus X$. Then for infinitely many $n$ we have that 
$a_n<b'$ and for every large $n$ that $b\le b_n$. We take a~$k$ such that $\frac{2}{k}\le 
b-b'$. Since $(a_n)\in C$, for every large $m$ and $n$ we have that $|a_m-a_n|\le\frac{b-b'}{2}$. Then for 
every large $n$ it holds that
$${\textstyle
a_n\le b'+\frac{b-b'}{2}=b-\frac{b-b'}{2}\le b_n-\frac{1}{k}\,.
}
$$
Hence $X<Y$.
\eproof

We show that
$\R\equiv (K,0_K,1_K,+,\cdot,<)$ is an ordered field. The ring structure follows 
easily from a~more general construction. Recall that $(R,0_R,1_R,+,\cdot)$ is 
a~commutative ring with the identity $1_R$, shortly 
a~{\em ring}, if 
the operations $+$ and $\cdot$ on the set $R$ are commutative and associative, the distinct elements
$0_R$ and $1_R$ of $R$ are neutral to $+$ and $\cdot$, respectively, the operation $\cdot$ is 
distributive to $+$, and every element $a\in R$ has a~unique additive inverse $-a\in R$ 
for which $a+(-a)=0_R$. The only missing axiom of a~{\em field} is that every 
$a\in R\setminus\{0_R\}$ has a~unique multiplicative 
inverse $\frac{1}{a}=1/a=a^{-1}\in R$ for which $a\cdot a^{-1}=1_R$. The two missing {\em 
order axioms} of an {\em ordered field} say that with a~linear order $(R,<)$, for every $a$, 
$b$ and $c$ in $R$ one has that
$$
a<b\Rightarrow a+c<b+c\,\text{ and }\,0_R<a\wedge 0_R<b\Rightarrow 0_R<a\cdot b\;.
$$
Recall that $U^{\N}$ denotes the set of sequences $(u_n)\sus U$. We give a~construction of a~ring as an $\N$-th power of another ring.

\begin{prop}\label{prop_NthPoRin}
Let 
$R_{\mathrm{ri}}\equiv (R,0_R,1_R,+,\cdot)$ 
be a~ring and $P\equiv R^{\N}$. Then the algebraic structure
$$
R_{\mathrm{ri}}^{\N}\equiv (P,\,0_P,\,1_P,\,+,\,\cdot)\,,
$$
where $0_P\equiv (0_R,0_R,\ds)$, $1_P\equiv (1_R,1_R,\ds)$ and the operations $+$ and $\cdot$ on $P$ are defined component-wise by the $+$ and $\cdot$ in $R_{\mathrm{ri}}$,  is a~ring. 
\end{prop}
\proof
The ring axioms hold in the ring $R_{\mathrm{ri}}^{\N}$ because they hold in every component, which is isomorphic to the ring $R_{\mathrm{ri}}$.
\eproof

\noindent
Note that $R_{\mathrm{ri}}^{\N}$ is never an integral domain, let alone a~field. 

We show that $\R$ is an ordered field. We build on the fact, 
taken for granted, that $\Q_{\mathrm{OF}}$ is an ordered field.

\begin{thm}\label{thm_KisOF}
The structure
$\R=(K,0_K,1_K,+,\cdot,<)$
is an ordered field.
\end{thm}
\proof
Let $\Q_{\mathrm{OF}}=(\Q,0,1,+,\cdot,<)$ be the ordered field of fractions, $Q\equiv \Q^{\N}$ and 
let
$$
\Q^{\N}_{\mathrm{ri}}\equiv (Q,\,0_Q,\,1_Q,\,+,\,\cdot)
$$
be the ring of Proposition~\ref{prop_NthPoRin}. 
Then $C\sus Q$ and we showed
above that $0_Q,1_Q\in C$ and that $C$ is closed to the operations $+$ and $\cdot$ on 
$Q$. Thus $C_{\mathrm{ri}}\equiv (C,0_Q,1_Q,+,\cdot)$ is a~subring
of $\Q^{\N}_{\mathrm{ri}}$. We have the surjection $G\cc C\to K$ such that
$G(a_n)=G(b_n)$ iff $(a_n)\sim(b_n)$. 
We obtained the neutral elements $0_K,1_K\in K$, the operations $+$ 
and $\cdot$ on $K$ and the relation $<$ on $K$ by pushing forward by $G$ 
the corresponding neutral elements, operations and relation in 
$C_{\mathrm{ri}}$. Satisfaction of the ring axioms is pushed forward 
too and therefore $\R$ is a ring.  It remains to prove: (i) the existence of 
multiplicative inverses, (ii) $(K,<)$ is a~linear order and (iii) the 
two order axioms.

(i) Let $X\in K\setminus\{0_K\}$, so that $X=G(a_n)$ with $(a_n)\in C$ and $(a_n)\not\sim 0_Q$. 
The latter means that for some $k$ 
and infinitely many $n$
one has $|a_n|\ge\frac{2}{k}$. Since $(a_n)$
is Cauchy, for every large $n$ it holds that  $|a_n|\ge\frac{1}{k}$. We define the sequence $(b_n)\in Q$ by
setting $b_n\equiv\frac{0}{1}$ if $a_n=\frac{0}{1}$, and $b_n\equiv \frac{1}{a_n}$ else. We show that $(b_n)\in C$. Since $(a_n)\in C$,
for any given $l$ we have for every large $m$ and $n$ that $|a_m-a_n|\le\frac{1}{l}$. Hence for every large $m$ and $n$,
$${\textstyle
|b_m-b_n|=\big|\frac{1}{a_m}-\frac{1}{a_n}\big|=\frac{|a_n-a_m|}{|a_m|\cdot|a_n|}\le\frac{k^2}{l}
}
$$
and $(b_n)$ is Cauchy. So $X^{-1}\equiv G(b_n)$ is a~cut and 
$X\cdot X^{-1}=G(a_nb_n)=G(\frac{1}{1},\frac{1}{1},\ds)=1_K$ because $(a_nb_n)\sim(\frac{1}{1},\frac{1}{1},\ds)$.

(ii) It follows from the definition of $<$  that $<$ is irreflexive and transitive. We show that $<$ is trichotomic. Suppose that $X\ne Y$ are cuts. Then there is an $a\in X\setminus Y$ or 
a~$b\in Y\setminus X$. We consider the former case, the latter is similar. It follows that in $(\Q,<)$ one has that $Y<a$.
Thus $Y\sus X$ and, by Proposition~\ref{prop_lesVerInc}, $X<Y$. 

(iii) Consider three cuts $X=G(a_n)$, $Y=G(b_n)$ and $Z=G(c_n)$, where $(a_n)$, $(b_n)$ and $(c_n)$ are in $C$. If
$X<Y$ then for some $k$ and every large $n$ it holds that $a_n\le b_n-\frac{1}{k}$. Thus for the same $n$, $a_n+c_n\le b_n+c_n-\frac{1}{k}$ and 
$X+Z=G(a_n+c_n)<G(b_n+c_n)=Y+Z$. If $0_K<X$ and $0_K<Y$ then for some 
$k$ and large $n$ one has that 
$\frac{1}{k}\le a_n$ and $\frac{1}{k}\le b_n$. So for the same $n$, 
$\frac{1}{k^2}\le a_nb_n$ and $0_K<G(a_nb_n)=X\cdot Y$.
\eproof

\noindent
It is easy to see that $\R$ is {\em Archimedean}, for every cut $X>0_K$ there 
is a~$k$ such that $X>E(\frac{1}{k})>0_K$.

\subsection{Real sequences}\label{subsec_sequences}

A~{\em real sequence $(x_n)$} is a~sequence in $K$, a~function from $\N$ to $K$. It is an HMC set. 
From now on we work in the linear order $(K,<)$. A~real number $x$ is 
an {\em upper bound of $(x_n)$} if  $x\ge x_n$ for every $n$. We say that $(x_n)$ is {\em bounded 
from above} if it has an upper bound. An upper bound $x$ of $(x_n)$ is the 
{\em least upper bound of $(x_n)$} if for every $y<x$ there is an $n$
such that $y<x_n$. We similarly define 
a~{\em lower bound of $(x_n)$}, {\em boundedness of $(x_n)$ from below} and 
the {\em largest lower bound of $(x_n)$}. The sequence $(x_n)$ is 
{\em bounded} if it is bounded both from above and below. We introduce 
suprema and infima of real sequences.

\begin{thm}\label{thm_supInR}
Any real sequence $(x_n)$ bounded from above has the least upper bound $x$. We call $x$ the 
supremum of $(x_n)$ and denote it by 
$\sup x_n$. Any real sequence $(x_n)$ bounded from below has the largest lower 
bound $y$. We call $y$ the infimum of $(x_n)$ and denote it by $\inf x_n$. If the elements $\sup x_n$ and $\inf x_n$ exist, they are unique.
\end{thm}
\proof
Uniqueness of suprema and infima follows from their definitions.
The replacement of $(x_n)$ with $(-x_n)$ shows that it suffices to 
prove just the first claim. Every $x_n$ is a~cut on $\Q$. We set
$${\textstyle
x\equiv\bigcup_{n=1}^{\infty}x_n\ \ (\sus\Q)
}
$$
and show that $x$ is a~cut and the least upper bound of $(x_n)$. 

We show that $x$ is a~cut. Clearly, $x\ne\emptyset$. 
If $y\in K$ is an upper bound of $(x_n)$ then, by 
Proposition~\ref{prop_lesVerInc}, $x_n\sus y$ for every $n$.
Thus also $x\sus y$ and $\Q\setminus x\ne\emptyset$ because 
$\Q\setminus y\ne\emptyset$. If $a<b$ are 
fractions and $b\in x$, then $b\in x_n$ for some $n$. Hence $a\in x_n$ and $a\in x$. Suppose that 
$a=\max(x)$ exists. Then $a\in x_n$ for some $n$, but then 
$a=\max(x_n)$ which is impossible. Thus $x$ has no largest element and 
is a~cut. 

By Proposition~\ref{prop_lesVerInc},  $x$ is an upper bound of $(x_n)$. 
Let $y$ be a~real number such that $y<x$. By 
Proposition~\ref{prop_lesVerInc}, $y\sus x$ and $y\ne x$. Thus there is a~fraction 
$a\in x\setminus y$ and $y<a$. Since $a\in x_n$ for some $n$, it holds that $y\sus x_n$ and $y\ne x_n$.
This means by Proposition~\ref{prop_lesVerInc} that $y<x_n$. So $x$ is the least 
upper bound of the sequence $(x_n)$.
\eproof

For a~real number $x$ its absolute value $|x|\equiv x$ if $x\ge0_K$, and $|x|\equiv-x$ if $x<0_K$.
The {\em triangle inequality} holds: for all real numbers $x_1$, $\ds$,
$x_n$,  $n\in\N$, it is true that
$|x_1+\ds+x_n|\le|x_1|+\ds+|x_n|$.
We remind limits of real sequences. 

\begin{defi}\label{def_realLimits}
A~real sequence $(x_n)$ has a~limit $x$, in symbols $\lim x_n=x$ or $\lim_{n\to\infty}x_n=x$, if for every $k$ there is an $n_0$ such that
$${\textstyle
n\ge n_0\Rightarrow|x_n-x|\le\frac{1}{k}\,.
}
$$
\end{defi} 
Recall that instead of $\frac{1}{k}$ one should write here more precisely 
$E(\frac{1}{k})$. Limits are unique. If a~real sequence $(x_n)$ 
has a~limit, we say that $(x_n)$ {\em converges}.

\begin{defi}\label{def_limInfin}
A~real sequence $(x_n)$ goes to $+\infty$ if for every $k$ 
there is an $n_0$ such that
$$
n\ge n_0\Rightarrow x_n\ge k\,.
$$
    
\end{defi}
Clearly, if $(x_n)$ goes to $+\infty$ then for every $y$ it holds that 
$\lim\frac{y}{x_n}=0_K$ (we ignore finitely many undefined fractions 
$\frac{y}{0}$). 
We prove four results on relations between limits and arithmetic in $\R$.

\begin{prop}\label{prop_limSeqDomi}
Suppose that $\lim y_n=0_K$ and that  for every large $n$ it holds that 
$|x_n|\le|y_n|$. Then $\lim x_n=0_K$.
\end{prop}
\proof
Let a~$k$ be given. Then for every large $n$,  
$|y_n|=|y_n-0_K|\le\frac{1}{k}$. So  for every large $n$ also $|x_n-0_K|=|x_n|\le |y_n|\le\frac{1}{k}$ and $\lim x_n=0_K$.
\eproof
\vspace{-3mm}
\begin{prop}\label{prop_linCombLimSeq}
Let $\lim x_n=x$ and $\lim y_n=y$. Then 
$\lim (z x_n+w y_n)=z x+w y$.
\end{prop}
\proof
Let a~$k$ be given. Then for every large $n$ one has that $|x-x_n|\le \frac{1}{k}$ 
and $|y-y_n|\le\frac{1}{k}$. The triangle inequality shows that for the same $n$,
$${\textstyle
|z x_n+w y_n-(z x+w y)|\le|z|\cdot|x_n-x|+|w|\cdot|y_n-y|
\le\frac{|z|+|w|}{k}\,.
}
$$
Hence $\lim (z x_n+w y_n)=z x+w y$.
\eproof
\vspace{-3mm}
\begin{prop}\label{prop_limProdu}
Let $\lim x_n=x$ and $\lim y_n=y$. Then $\lim x_ny_n=xy$.
\end{prop}
\proof
Convergent sequences are bounded and therefore for some $z>0_K$ and every $n$ it holds that that $|y_n|\le z$. For a~given $k$ it holds for every large $n$ that
$|x-x_n|\le\frac{1}{k}$ and $|y-y_n|\le\frac{1}{k}$. By the triangle inequality it holds for the same $n$ that 
$${\textstyle
|x_ny_n-xy|\le|x_n-x|\cdot|y_n|+|x|\cdot|y_n-y|\le \frac{z}{k}+\frac{|x|}{k}=\frac{z+|x|}{k}\,.
}
$$
Hence $\lim x_ny_n=xy$.
\eproof
\vspace{-3mm}
\begin{prop}\label{prop_OrderAndLimSeq}
Let $\lim x_n=x$, $\lim y_n=y$ and let for every $n\ge n_0$ it hold that  
$x_n\le y_n$. Then $x\le y$.
\end{prop}
\proof
It follows that for every $k$ we have for large $n$ that $x-\frac{1}{k}\le x_n\le y_n\le y+\frac{1}{k}$. Sending $k\to\infty$ we get that $x\le y$.
\eproof

If $x_1\le x_2\le\ds$, the sequence {\em $(x_n)$ increases}, and if $x_1\ge 
x_2\ge\ds$ then it {\em decreases}. When the inequalities are strict  then $(x_n)$ {\em strictly 
increases}, resp. {\em strictly decreases}. We say that
$(x_n)$ is {\em monotone} if it increases or 
decreases. We show that every bounded monotone sequence converges.

\begin{prop}\label{prop_monoSeq}
If $(x_n)$ increases and is bounded from above, resp. decreases and is bounded from
below, then it converges and $\lim x_n=\sup x_n$, resp. $\lim x_n=\inf x_n$.
\end{prop}
\proof
We assume that $(x_n)$ increases and is bounded from above, the other case is treated similarly. Using Theorem~\ref{thm_supInR} we set 
$x\equiv\sup x_n$. Let a~$k$ be given. By the definition of supremum there is an $l$ such that $x_l>x-\frac{1}{k}$. Then for any $n\ge l$,
$${\textstyle
x-\frac{1}{k}<x_l\le x_n\le x<x+\frac{1}{k}\,\text{ and }\,|x_n-x|\le\frac{1}{k}\,.
}
$$
Hence $\lim x_n=x$.
\eproof

We determine limits of exponential sequences.

\begin{cor}\label{cor_limXnaN}
The limit $\lim x^n$ is $0_K$ for $|x|<1_K$, $1_K$ for $x=1_K$ and does not exist for $x\le-1_K$ or $x>1_K$.
\end{cor}
\proof
Let $|x|<1_K$. Since $|x^n|=|x|^n$ and $\lim 0_K^n=0_K$, we may assume that $0_K<x<1_K$. Then $(x^n)$ decreases and is bounded from below by $0_K$. Using the two previous propositions we get that
$\lim x^n=\inf x^n\equiv y\ge0_K$. 
Suppose for the contrary that $y>0_K$. By the definition of infimum there is an $n$ such that 
$x^n<\frac{y}{x}$. But then $x^{n+1}<y$, contradicting that $y$ is a~lower bound of $(x^n)$. Thus 
$y=0_K$.

For $x=1_K$ the limit is trivial. For $x>1_K$ we show similarly that $(x^n)$ is not bounded from 
above and so $\lim x^n$ does not exist. If $x\le-1_K$ then $x^n\le-1_K$ for odd $n$, $x^n\ge1_K$ for
even $n$ and again $\lim x^n$ does 
not exist.
\eproof

\noindent
For $x>1_K$ the sequence $(x^n)$ goes to $+\infty$.

Factorial growth trumps every exponential one.

\begin{cor}\label{cor_limSeqExpFac}
For any real numbers $y$ and $z$, 
$\lim_{n\to\infty}\frac{y z^n}{n!}=0_K$. 
\end{cor}
\proof
It suffices to prove that for every $z>0_K$ we have  
$\lim_{n\to\infty}\frac{z^n}{n!}=0_K$. We take an $m\in\N$ such that $m>z$ and  consider the sequence
$(x_n)\equiv(\frac{z^m}{m!}, \frac{z^{m+1}}{(m+1)!},\ds)$. It has positive terms and decreases. By Propositions~\ref{prop_OrderAndLimSeq} and \ref{prop_monoSeq},  
$${\textstyle
\lim_{n\to\infty}\frac{z^n}{n!}=\lim x_n=\inf x_n\equiv y'\ge0_K\,.
}
$$
Suppose for the contrary that $y'>0_K$. By the definition of infimum, there is an $n$ 
such that $0_K<x_n<\frac{y'm}{z}$. But then $x_{n+1}=x_n\frac{z}{m+n}<y'$, contradicting 
that $y'$ is a~lower bound of $(x_n)$. Thus $y'=0_K$. 
\eproof

We prove the infinite form of the {\em Erd\H os--Szekeres lemma}. Recall that $(y_n)$ 
is a~{\em subsequence} of 
a~sequence $(x_n)$ if for some strictly increasing sequence $(m_n)\sus\N$ it 
holds for every $n$ that $y_n=x_{m_n}$.

\begin{lem}\label{lem_mono}
Every real sequence $(x_n)$ has a~monotone 
subsequence.
\end{lem}
\proof
Let $(x_n)$ be a~real sequence and $H\equiv\{m\in\N:\;m<n \Rightarrow x_m<x_n\}$. If $H$ is infinite, $H=\{m_1<m_2<\ds\}$, then 
$(x_{m_n})$ is a~strictly increasing subsequence 
of $(x_n)$. If $H$ is finite, we take any $m_1>\max(H)$ (or any $m_1$ if $H=\emptyset$). Since 
$m_1\not\in H$, there is an $m_2$ such that $m_1<m_2$ and $x_{m_1}\ge x_{m_2}$. Since  
$m_2\not\in H$, there is an $m_3$ such that $m_2<m_3$ and $x_{m_2}\ge x_{m_3}$, and so 
on. Hence $(x_{m_n})$ is a~decreasing subsequence of $(x_n)$.
\eproof

A~real sequence $(x_n)$ is {\em Cauchy} if for every $k$ there is an $n_0$ such that
$${\textstyle
m,\,n\ge n_0\Rightarrow|x_m-x_n|\le \frac{1}{k}\,.
}
$$
We prove the {\em metric completeness} of $\R$. 

\begin{thm}\label{thm_complR}
A~real sequence $(x_n)$ is Cauchy if and only if it converges.
\end{thm}
\proof
Let $(x_n)$ be Cauchy. Thus $(x_n)$ is bounded and by the previous lemma has 
a~(bounded) monotone subsequence $(x_{m_n})$. Using Proposition~\ref{prop_monoSeq} we take $y\equiv\lim x_{m_n}$ and show that $\lim x_n=y$. Let a~$k$ be given. Then we have for every large $m$ and $n$ that $|x_{m_n}-y|\le\frac{1}{2k}$ and $|x_m-x_n|\le\frac{1}{2k}$. Then (since 
$m_n\ge n$) one has for the same $n$ that 
$${\textstyle
|x_n-y|\le|x_n-x_{m_n}|+|x_{m_n}-y|\le\frac{1}{2k}+\frac{1}{2k}=\frac{1}{k}\,.
}
$$
Hence $\lim x_n=y$.

Let $\lim x_n=x$ and a~$k$ be given. Then it holds for every large $n$ that $|x_n-x|\le\frac{1}{2k}$. We have by the triangle inequality for the same $m$ and $n$ that 
$${\textstyle
|x_m-x_n|\le|x_m-x|+|x-x_n|\le \frac{1}{2k}+\frac{1}{2k}=\frac{1}{k}\,.
}
$$
Hence $(x_n)$ is Cauchy.
\eproof

We conclude this subsection with the {\em Bolzano--Weierstrass theorem}.

\begin{thm}\label{thm_BW}
Any bounded real sequence has a~convergent subsequence.
\end{thm}
\proof
This follows from Lemma~\ref{lem_mono} and  Proposition~\ref{prop_monoSeq}.
\eproof

\subsection{Real series}\label{subsec_series}

We define series and their sums.

\begin{defi}\label{def_sumSer}
A~series is a~real sequence $(x_0,x_1,\ds)$, a~function from 
$\omega$ to $K$. One denotes it also as $\sum_{n=0}^{\infty}x_n$.
It converges if the limit
$$
\lim_{n\to\infty}(x_0+x_1+\ds+x_n)
$$
exists. We call this limit the sum of the series and denote it again 
by $\sum_{n=0}^{\infty}x_n$.
Let $m\in\N$ and $x_m$, $x_{m+1}$, $\ds$ be real numbers. 
We denote by 
$\sum_{n=m}^{\infty}x_n$ the series $\sum_{n=0}^{\infty}y_n$ defined by 
$y_0=\ds=y_{m-1}\equiv 0_K$ and $y_n\equiv x_n$ for $n\ge m$.
\end{defi}
Every series is an HMC set. We review absolutely convergent series.

\begin{defi}\label{def_absCon}
A~series $\sum_{n=0}^{\infty}x_n$ absolutely converges if the series
$\sum_{n=0}^{\infty}|x_n|$
converges.    
\end{defi}

Linear combinations of series preserve sums.

\begin{prop}\label{prop_linSums}
Let we have sums $\sum_{n=0}^{\infty}x_n=x$ and 
$\sum_{n=0}^{\infty}y_n=y$. Then we have the sum
$\sum_{n=0}^{\infty}(z x_n+w y_n)=zx+wy$.
\end{prop}
\proof
This follows from Definition~\ref{def_sumSer} and Proposition~\ref{prop_linCombLimSeq}.
\eproof

Absolute convergence implies convergence.

\begin{prop}\label{prop_ACimplC}
If a~series $\sum_{n=0}^{\infty}x_n$ absolutely converges then it converges, and for the sums the infinite triangle inequality
$\big|\sum_{n=0}^{\infty}x_n\big|\le\sum_{n=0}^{\infty}|x_n|$ holds.
\end{prop}
\proof
Suppose that the series $\sum_{n=0}^{\infty}x_n$ absolutely converges and let $s_n\equiv\sum_{j=0}^n x_n$ and $t_n\equiv\sum_{j=0}^n|x_n|$. By Theorem~\ref{thm_complR} 
the sequence $(t_n)$ is Cauchy: for a~given $k$ for every large numbers 
$m\le n$ it holds that
$${\textstyle
\big||x_m|+|x_{m+1}|+\ds+|x_n|\big|=|x_m|+|x_{m+1}|+\ds+|x_n|\le\frac{1}{k}\,.
}
$$
By the triangle inequality for the same numbers $m\le n$ it holds that
$${\textstyle
|x_m+x_{m+1}+\ds+x_n|\le|x_m|+|x_{m+1}|+\ds+|x_n|\le\frac{1}{k}\,.
}
$$
Hence $(s_n)$ is Cauchy.  By Theorem~\ref{thm_complR} the series $\sum_{n=0}^{\infty}x_n$ converges.

For every $n$ one has that $|s_n|\le t_n$, equivalently $-t_n\le s_n\le t_n$. By Propositions~\ref{prop_linCombLimSeq} and \ref{prop_OrderAndLimSeq} and by the limit transition $n\to\infty$ we get for sums
the same inequalities 
$${\textstyle
-\sum_{n=0}^{\infty}|x_n|\le
\sum_{n=0}^{\infty}x_n\le 
\sum_{n=0}^{\infty}|x_n|\,.
}
$$
Hence $|\sum_{n=0}^{\infty}x_n|\le\sum_{n=0}^{\infty}|x_n|$.
\eproof

We sum the {\em geometric series} $\sum_{n=m}^{\infty}x^n$.

\begin{prop}\label{prop_geomSer}
Let $m\in\omega$. The series $\sum_{n=m}^{\infty}x^n$ converges iff $|x|<1_K$ and then it has the sum 
$\sum_{n=m}^{\infty}x^n=\frac{x^m}{1_K-x}$,
\end{prop}
\proof
Let $x\ne1_K$ and $n\ge m$. Then $s_n\equiv\sum_{j=m}^n x^j=x^m\cdot\frac{1_K-x^{n-m+1}}{1_K-x}$. By 
Proposition~\ref{prop_linCombLimSeq} and Corollary~\ref{cor_limXnaN}, if $|x|<1_K$ 
then $\lim s_n=\frac{x^m}{1_K-x}$, and if $|x|\ge 1_K$ then the limit does not exist. For 
$x=1_K$ one has that $s_n=n-m+1_K$ 
and $\lim s_n$ does not exist.
\eproof

We have the following simple criterion of convergence of series.

\begin{prop}\label{prop_dominSeries}
Suppose that $x$ is in $K$ and that $\sum_{n=0}^{\infty}x_n$ and
$\sum_{n=0}^{\infty}y_n$ are series such that $\sum_{n=0}^{\infty}y_n$ converges and $|x_n|\le xy_n$ for every $n\ge m$. Then $\sum_{n=0}^{\infty}x_n$ converges. 
\end{prop}
\proof
For every $n$ with $n\ge m$,
$$
|x_0|+\ds+|x_n|\le|x_0|+\ds+|x_{m-1}|+x\cdot(y_0+\ds+y_n)\;.
$$
So $(|x_0|+\ds+|x_n|)$ increases and is bounded from 
above; it converges by Proposition~\ref{prop_monoSeq}. 
Hence $\sum_{n=0}^{\infty}x_n$ absolutely converges, and converges 
by Proposition~\ref{prop_ACimplC}. 
\eproof

We review the Cauchy products of series.

\begin{thm}\label{thm_prodACser}
Let 
$S=\sum_{n=0}^{\infty}x_n$ and $T=\sum_{n=0}^{\infty}y_n$ be absolutely convergent series
with respective sums $x$ and $y$. Then their Cauchy product
$${\textstyle
S\odot T=\sum_{n=0}^{\infty}z_n\equiv\sum_{n=0}^{\infty}(x_0y_n+x_1y_{n-1}+\ds+x_ny_0)
}
$$
has the sum $x y$ and absolutely converges.
\end{thm}
\proof
Let $s_n\equiv\sum_{j=0}^n x_j$, $t_n\equiv\sum_{j=0}^n y_j$ and $u_n\equiv\sum_{j=0}^n z_j$. By the triangle inequality, 
$$
|u_n-xy|\le |u_n-s_nt_n|+|s_nt_n-xy|\equiv A_n+B_n\;.
$$ 
By Proposition~\ref{prop_limProdu}, $\lim B_n=0_K$. We show that also $\lim A_n=0_K$. 
We consider sums $x'\equiv\sum_{n=0}^{\infty}|x_n|$ and 
$y'\equiv\sum_{n=0}^{\infty}|y_n|$. Since $u_n=\sum_{i+j\le n}x_iy_j$ (with $i,j,n\in\omega$), 
we get that $A_n=\big|-\sum_{\substack{i,\,j\le n\\ i+j>n}}
x_iy_j\big|$. Thus
$${\textstyle
A_n\le\big(\sum_{n/2<i\le n}|x_i|\big)\cdot y'+x'\cdot\big(\sum_{n/2<j\le n}|y_j|\big)\equiv C_n\cdot y'+x'\cdot D_n\;.
}
$$
By the absolute convergence of $S$ and $T$ the sequences
$(|x_0|+\ds+|x_n|)$ and $(|y_0|+\ds+|y_n|)$
are Cauchy. Hence $\lim C_n=\lim D_n=0_K$. By Proposition~\ref{prop_linCombLimSeq}, 
$\lim A_n=0_K$ and $\lim u_n=xy$. We have shown that $S\odot T$ converges. This 
proof shows also that the sum
$$
{\textstyle
\sum_{n=0}^{\infty}|z_n|\le 
\sum_{n=0}^{\infty}(|x_0|\cdot|y_n|+
|x_1|\cdot|y_{n-1}|+\ds+|x_n|\cdot|y_0|)=x'\cdot y'\;.
}
$$ 
Thus $S\odot T$ absolutely converges. 
\eproof
\vspace{-3mm}
\subsection{Countable functions}\label{subsec_HMCfun}

Now we swerve from the main road
of real analysis and start to develop a~theory of countable univariate real 
functions. For any set $M\sus\Q$ we define
$$
\mathcal{F}(M)\equiv\{f:\;f\cc M\to K\}\,.
$$
If we 
write $\mathcal{F}(M)$, $\mathcal{F}(N)$, $\ds$, the sets $M$, $N$, 
$\ds$ are always subsets of 
$\Q$. Every function $f\in\mathcal{F}(M)$ is at most countable, equivalently HMC. We define how functions are added, multiplied and divided.

\begin{defi}\label{def_aritOfFun}
Let $f\in\mathcal{F}(M)$ and $g\in\mathcal{F}(N)$. Then $f+g,fg\cc M\cap N\to K$ have values $(f+g)(a)\equiv 
f(a)+g(a)$ and $(fg)(a)\equiv f(a)g(a)$. As for division, $f/g\cc M\cap N\setminus Z(g)\to K$ has values 
$(f/g)(a)\equiv\frac{f(a)}{g(a)}$, where $Z(g)\equiv\{a\in N:\;g(a)=0_K\}$.
\end{defi}
In this article we do not use division of functions.

For any $x\in K$ we define $k_x\in\mathcal{F}(\Q)$ to be the 
{\em constant function} with the  values $k_x(a)=x$. By 
$\mathrm{id}\in\mathcal{F}(\Q)$ we denote the {\em identity function}
$\mathrm{id}(a)=a$. We simplify notation by writing $x$
instead of $k_x$, and $a$ instead of $\mathrm{id}(a)$. For 
example, in Proposition~\ref{prop_derExp}, 
$c\cdot\mathrm{e}^{ca}$ means $k_c(a)\cdot\mathrm{e}^{k_c(a)\cdot
\mathrm{id}(a)}$. We define polynomials.

\begin{defi}\label{def_polyn}
A~polynomial is any function $p\in\mathcal{F}(\Q)$ for which there exists an $n$-tuple, $n\in\N$, of 
functions $f_1,\ds,f_n$ in $\mathcal{F}(\Q)$ such that $p=f_n$ and for every $i\in[n]$ the 
function $f_i$ is a~constant function or the identity function or for some $j,k\in[i-1]$ it holds that $f_i=f_j+f_k$
or $f_i=f_jf_k$.
\end{defi}
By a~polynomial in $\mathcal{F}(M)$
we mean the restriction of a~polynomial to $M$. It is not hard to show that every polynomial $p$ is either the {\em zero polynomial} $p=k_{0_K}$ or has the unique {\em canonical form}
$$
{\textstyle
p=p(a)=\sum_{j=0}^n x_ja^j=\sum_{j=0}^n k_{x_j}(a)\cdot\mathrm{id}(a)^j\,,
}
$$
where $a$ runs in $\Q$, $n\in\omega$, $x_j\in K$ and $x_n\ne0_K$. 
It is more standard to define polynomials by these canonical forms. We set $\deg k_0\equiv-\infty$ and if 
a~polynomial $p$ has the above canonical form then the {\em degree} of $p$ 
is $\deg p\equiv n$.

In our approach to real analysis uniform continuity is fundamental.

\begin{defi}\label{def_unifCont}
A~function $f\in\mathcal{F}(M)$ is uniformly continuous (on $M$), abbreviated {\em 
UC}, if for every $k$ there is an $l$ such that for any $a,b\in M$, 
$${\textstyle
|a-b|\le \frac{1}{l}\Rightarrow |f(a)-f(b)|\le\frac{1}{k}\,.
}
$$
\end{defi}
It is easy to see that constant functions $k_x(a)$ and the identity 
function $\mathrm{id}(a)$ are UC on any $M\sus\Q$.
If $f\cc M\to K$ with $M\sus K$, which is allowed only in 
Theorem~\ref{thm_extending}, then uniform continuity of $f$ is 
defined in the same way.

For any function 
$f\in\mathcal{F}(M)$ we define
its absolute value $|f|\in\mathcal{F}(M)$ by 
$|f|(a)\equiv|f(a)|$. This construction 
preserves uniform continuity.

\begin{prop}\label{prop_UCandAbsVal}
If $f\in\mathcal{F}(M)$ is {\em UC} then so is $|f|$ ($\in\mathcal{F}(M)$).
\end{prop}
\proof
Let $f$ be as stated and let a~$k$ be given. We take an $l$ such 
that for any $a,b\in M$ it holds that 
$|a-b|\le \frac{1}{l}\Rightarrow|f(a)-f(b)|\le
\frac{1}{k}$. Suppose that $a,b\in M$ satisfy $|a-b|\le\frac{1}{l}$. If $f(a)$ and $f(b)$ have equal signs, 
or one of them is $0_K$, then
$${\textstyle
\big||f|(a)-|f|(b)\big|=|(\pm1_K)\cdot(f(a)-f(b))|=|f(a)-f(b)|\le\frac{1}{k}\,.
}
$$
If they have different signs then
$|f(a)-f(b)|\le\frac{1}{k}$ implies that $|f(a)|,|f(b)|\le\frac{1}{k}$.
Thus again
$${\textstyle
\big||f|(a)-|f|(b)\big|\le\max(\{|f(a)|,\,|f(b)|\})\le\frac{1}{k}\,.
}
$$
We see that $|f|$ is UC.
\eproof

A~set $U\sus K$ is {\em bounded} if for some $y\ge0_K$ and all $u\in U$ we have that $|u|\le y$; then we say that 
{\em $y$ bounds $U$}. A~function $f\in\mathcal{F}(M)$ is {\em bounded} if its image $f[M]$ is 
bounded. We show that any UC function with bounded definition domain $M\sus\Q$ is bounded.

\begin{prop}\label{prop_UCbound}
If $f\in\mathcal{F}(M)$ is {\em UC} and $M$ is bounded then $f$ is bounded.
\end{prop}
\proof
By the assumption there is an $l$ such that for any $a,b\in M$ it holds that
$|a-b|\le\frac{1}{l} \Rightarrow|f(a)-f(b)|\le 1_K$. Since $M$ is bounded, there exists a~finite set $N\sus M$ such that for every $a\in 
M$ there is a~$b_a\in N$ with $|a-b_a|\le\frac{1}{l}$. Let $y$ bound the finite set 
$f[N]$. Then for any $a\in M$,
$$
|f(a)|\le|f(a)-f(b_a)|+|f(b_a)|\le 1_K+y\,.
$$
Hence $f$ is bounded.
\eproof

Linear combinations and products preserve uniform continuity.

\begin{prop}\label{prop_UCandSums}
If $f,g\in\mathcal{F}(M)$ are {\em UC} then $xf+yg$ ($\in\mathcal{F}(M)$) is {\em UC}.
\end{prop}
\proof
Let $f$ and $g$ be as stated and $x,y\in K$. For a~given $k$ we take an $l$ such that for any 
$a,b\in M$ it holds that $|a-b|\le \frac{1}{l}$ $\Rightarrow$ $|f(a)-f(b)|,|g(a)-g(b)|\le\frac{1}{k}$. 
Then for every $a,b\in M$ with $|a-b|\le \frac{1}{l}$ we have that
$${\textstyle
|(xf+yg)(a)-(xf+yg)(b)|\le 
|x|\cdot|f(a)-f(b)|+|y|\cdot|g(a)-g(b)|\le\frac{|x|+|y|}{k}\,.
}
$$
Hence $xf+yg$ is UC.
\eproof
\vspace{-3mm}
\begin{prop}\label{prop_UCandPr}
If $f,g\in\mathcal{F}(M)$ are {\em UC} and $M$ is bounded, then the product $fg$ ($\in \mathcal{F}(M)$) is {\em UC}.
\end{prop}
\proof
Let $f$, $g$ and $M$ be as stated,  and let a~$k$ be given. Using 
Proposition~\ref{prop_UCbound} we take a~constant $y$ that bounds 
both $f[M]$ and $g[M]$. We take an $l$ such that for any $a,b\in M$ it holds that $|a-b|\le\frac{1}{l}$ $\Rightarrow$ 
$|f(a)-f(b)|,|g(a)-g(b)|\le \frac{1}{k}$. Then for every $a,b\in M$ 
with $|a-b|\le \frac{1}{l}$ we have that
$${\textstyle
|(fg)(a)-(fg)(b)|\le|f(a)|\cdot|g(a)-g(b)|+|g(b)|\cdot|f(a)-f(b)|\le\frac{2y}{k}\,.
}
$$
Hence $fg$ is UC.
\eproof

Using Definition~\ref{def_polyn} and the last two propositions, we get the next corollary.

\begin{cor}\label{cor_XnaNisUC}
If $M\sus\Q$ is bounded, then
any polynomial $p\in\mathcal{F}(M)$ is {\em UC}.
\end{cor}

Recall that if $M\sus\Q$, then $x\in K$
lies in the closure $\overline{M}$ (i.e., $\overline{E[M]}$) of $M$
iff $x=\lim x_n$ for some $(x_n)\sus E[M]$. In our approach to real analysis the next theorem is fundamental.

\begin{thm}\label{thm_extending}
Every {\em UC} function $f\in
\mathcal{F}(M)$ with $M\sus\Q$ has a~unique continuous extension $g\cc M\cup\{x\}\to K$ to every point $x\in\overline{M}$. 
\end{thm}
\proof
Let $f$, $M$ and $x$ be as stated. The uniqueness of the continuous extension of 
$f$, if it exists, follows from its continuity. We show that it exists. Let 
$(a_n)\sus M$ be such that $\lim a_n=x$ and a~$k$ be given. We take an $l$ such 
that for any $a,b\in M$ with $|a-b|\le\frac{1}{l}$ one has that 
$|f(a)-f(b)|\le 
\frac{1}{k}$. By Theorem~\ref{thm_complR} the sequence $(a_n)$ is Cauchy. Thus  for 
every large $m$ and $n$ it holds that $|a_m-a_n|\le \frac{1}{l}$. For the same 
$m$ and $n$ we have that $|f(a_m)-
f(a_n)|\le \frac{1}{k}$. Thus $(f(a_n))$ is Cauchy. By Theorem~\ref{thm_complR}, 
$\lim f(a_n)=y$ for some $y$. We set $g(x)\equiv y$. Let $(b_n)\sus M$ be 
another sequence with $\lim b_n=x$. Then the sequence
$(c_n)\equiv (a_1,b_1,a_2,b_2,\ds)$ has the limit $x$. Since $(f(c_n))$ 
contains $(f(a_n))$ and $(f(b_n))$ as subsequences,
$\lim f(a_n)=\lim f(c_n)=\lim f(b_n)$
and $\lim f(a_n)=\lim f(b_n)=y$. Thus the value $g(x)$ is correctly defined as it 
does not depend on the sequence $(a_n)$ and, by Heine's definition of pointwise continuity, $g$ is continuous at $x$.
\eproof

\noindent
We call the argument proving the uniqueness of the limit $y$ the {\em 
interleaving argument};
later we use it again. Abusing notation, from now on we denote the extending 
function $g$ as $f$. It appears that it would be more elegant to consider the whole UC 
extension $f\cc\overline{M}\to K$, but this function is in general uncountable.

We introduce functional composition of countable functions.

\begin{defi}\label{def_composFunct}
Let $M,N\sus\Q$, $g\in\mathcal{F}(N)$ and let $f\in\mathcal{F}(M)$ be {\em UC}. We define the composite function $f(g)\in\mathcal{F}(P)$ by setting $P\equiv\{a\in N:\;g(a)\in\overline{M}\}$ and $f(g)(a)\equiv f(g(a))$, where $f(g(a))$ is the value at $g(a)$ of the extension of $f$ provided by Theorem~\ref{thm_extending}.
\end{defi}
It is not hard to show that if also $g$ is UC then $f(g)$ is UC. 

The min-max principle for countable functions takes a~simple form.

\begin{thm}\label{thm_attaGlobExtr}
Suppose that $M\sus\Q$ is bounded
and $f\in\mathcal{F}(M)$ is {\em UC}. Then there exist 
$y,y'\in\overline{M}$ such that $f(y)\le f(a)\le f(y')$ for every 
$a\in M$.
\end{thm}
\proof
Let $M$ and $f$ be as stated. We show that $y$ exists, the existence of $y'$ is proven 
similarly. We take any sequence $(a_n)\sus\Q$ such that $M=\{a_n:\;n\in\N\}$.
By Proposition~\ref{prop_UCbound} the sequence $(f(a_n))$ is bounded. Using Theorem~\ref{thm_supInR} we set $z\equiv \inf f(a_n)$. 
By the definition of infimum there is a~sequence $(b_n)\sus M$ such that $z=\lim f(b_n)$. 
By Theorem~\ref{thm_BW} there exist $y\in K$ and a~subsequence $(b_{m_n})$ of $(b_n)$ such that 
$\lim_{n\to\infty} b_{m_n}=y$. 
Hence $y\in\overline{M}$. By Theorem~\ref{thm_extending} the extended value $f(y)=z$  
because $(f(b_{m_n}))$ is a~subsequence of $(f(b_n))$. Since $z$ is a~lower bound of the 
sequence $(f(a_n))$, we get that $f(y)=z\le f(a)$ for every $a\in M$. 
\eproof

\noindent
Clearly, $f(y)\le f(x)\le f(y')$ even for every $x\in\overline{M}$.

We do not need the next countable version of the {\em 
Bolzano--Cauchy theorem} and its proof 
is left as an exercise for the reader. It shows that in countable 
real analysis the total disconnectedness of nontrivial 
rational intervals is not a~problem.

\begin{thm}\label{thm_IVT}
Let $u<v$ be real numbers and $I\equiv [u,v]_{\Q}$. If $f\in\mathcal{F}(I)$ is {\em UC} 
and $f(u)f(v)\le0_K$, then there is an $x\in K$ such that $u\le x\le v$
and $f(x)=0_K$.
\end{thm}

\subsection{Derivatives of countable functions}

In this key section, we obtain countable versions of classical 
results on derivatives of univariate real functions. 
We achieve it with the help of uniform derivatives.

Derivatives of functions in $\mathcal{F}(M)$ are defined as for standard real functions. A~point $a\in M$ for $M\sus\Q$ is an {\em 
isolated point of $M$} if there is a~$k$ such that $(a-\frac{1}{k},a+
\frac{1}{k})_{\Q}\cap M=\{a\}$. Else $a$ is a~{\em limit point of $M$}. Clearly, $a$ 
is a~limit point of $M$ iff $a=\lim a_n$ for some $(a_n)\sus M\setminus\{a\}$. We denote the set of limit points of $M$ ($\sus\Q$) by $L(M)$ ($\sus M$).

\begin{defi}\label{def_deriv}
Let $M\sus\Q$, $f\in\mathcal{F}(M)$ and $a\in L(M)$. If $y\in K$ has the property that for every $k$ there is an $l$ such that for any $b\in M$ the implication
$${\textstyle
0_K<|a-b|\le\frac{1}{l}\Rightarrow
\big|\frac{f(b)-f(a)}{b-a}-y\big|\le\frac{1}{k}
}
$$
holds, we define $f'(a)\equiv y$ and say that $f'(a)$ is the derivative of $f$ at $a$.
\end{defi}
We allow only values in $K$ for derivatives
(not $\pm\infty$). For $f\in\mathcal{F}(M)$ we define the set
$$
D(f)\equiv \{a\in L(M)\;|\;\text{$f'(a)$ exists}\}\ \ (\sus L(M))\,.
$$
It is the definition domain of the function $f'\cc D(f)\to K$
that sends any $a$ to $f'(a)$. We call $f'$ the {\em 
derivative of $f$}. No isolated point of $M$ 
lies in $D(f)$. As for classical functions, the restriction $f\,|\,D(f)$ 
is continuous. But now this does not suffice because we need
UC functions. For example, consider 
the function $f\in\mathcal{F}([0,1]_{\Q})$ given as $f(a)=0_K$
for $a<\frac{1}{\sqrt{2}}$, and $f(a)=1_K$ for 
$a>\frac{1}{\sqrt{2}}$. Then $D(f)=[0,1]_{\Q}$, $f'=0_K$ and $f$ is continuous, but $f$ is not UC. To avoid such situations we introduce {\em 
uniform derivatives}.

\begin{defi}\label{def_unifDeri}
Let $M\sus\Q$, $f\in\mathcal{F}(M)$ and $N\sus D(f)$. We say that the derivative $f'\cc D(f)\to 
K$ is uniform on $N$ if for every $k$ there is an $l$ such that for every $a\in N$ and every $b\in M$ with $0_K<|a-b|\le\frac{1}{l}$ it holds that
$${\textstyle
\big|\frac{f(b)-f(a)}{b-a}-f'(a)\big|\le\frac{1}{k}\,.
}
$$
For $N=D(f)$ we simply say that $f'$ is uniform. 
\end{defi}
For $b\ne a\in D(f)$ we define $\de_{a,f}(b)\equiv
\frac{f(b)-f(a)}{b-a}
-f'(a)$, and for $b=a\in D(f)$ we set 
$\de_{a,f}(a)\equiv0_K$. 

We compute derivatives of constant functions and the identity.

\begin{prop}\label{prop_UderKonIde}
Let $M\sus\Q$. For any $x\in K$ the derivative $(k_x\,|\,M)'=k_{0_K}\,|\,L(M)$ and is uniform. The 
derivative $(\mathrm{id}\,|\,M)'=k_{1_K}\,|\,L(M)$ and is uniform.     
\end{prop}
\proof
Let $x\in K$ and $M\sus\Q$.
It is easy to compute that $(k_x\,|\,M)'=k_{0_K}\,|\,L(M)$ and that $(\mathrm{id}\,|\,M)'=k_{1_K}
\,|\,L(M)$. Let a~$k$ be given. We have for any $a\in L(M)$ and any $b\in M$ with $0_K<|a-b|\le1_K$ that
$$
{\textstyle
\big|\frac{k_x(b)-k_x(a)}{b-a}-k_{0_K}(a)\big|=
\big|\frac{x-x}{b-a}-0_K\big|
=0_K\le\frac{1}{k}\,,
}
$$
and that
$$
{\textstyle
\big|\frac{\mathrm{id}(b)-\mathrm{id}(a)}{b-a}-k_{1_K}(a)\big|=\big|\frac{b-a}{b-a}-1_K\big|=0_K\le\frac{1}{k}\,.
}
$$
Hence all these derivatives are uniform.
\eproof
 
Uniformity of $f'$ implies, under mild restrictions, that 
$f\,|\,D(f)$ is UC.

\begin{prop}\label{prop_UCandDer}
Let $M\sus\Q$. For any function $f\in \mathcal{F}(M)$ with uniform derivative $f'$ ($\in\mathcal{F}(D(f))$) 
the following hold. 
\begin{enumerate}
    \item If $f'$ is bounded then $f\,|\,D(f)$ is {\em UC}.
    \item If $M$ is bounded and $f'$ is {\em UC} then $f\,|\,D(f)$ is {\em UC}.
\end{enumerate}
\end{prop}
\proof
Let $M$ and $f$ be as stated. 1. Let $y$ bound $f'[D(f)]$ and a~$k$ 
be given. We take an $l$, $l\ge k$, such that the condition in 
Definition~\ref{def_unifDeri} holds, with $1_K$ in place of 
$\frac{1}{k}$. Then for any $a,b\in D(f)$ with $|a-b|\le\frac{1}{l}$,
$${\textstyle
|f(b)-f(a)|\le\big(|f'(a)|+|\de_{a,f}(b)|\big)\cdot|a-b|\le(y+1_K)\cdot|a-b|\le
\frac{y+1_K}{k}\,.
}
$$
Hence $f$ is UC.

2. If $M$ is bounded and $f'$ is UC then by Proposition~\ref{prop_UCbound} $f'$ is bounded and we 
can use part~1.
\eproof

We differentiate linear combinations.

\begin{prop}\label{prop_linDer}
Let $f,g\in\mathcal{F}(M)$ have uniform derivatives and $x,y\in K$. Then
$D(f)\cap D(g)\sus D(xf+yg)\sus M$, the restricted derivative
$$
(xf+yg)'\,|\,D(f)\cap D(g)=xf'+yg'
$$
and is uniform on $D(f)\cap D(g)$. 
\end{prop}
\proof
Let $f$, $g$, $M$, $x$ and $y$ be as stated and let a~$k$ be given.
We take an $l$ such that for any $a\in D(f)$, resp. $a\in D(g)$, and any $b\in M$ with $0_K<|a-b|\le\frac{1}{l}$ it holds that
$${\textstyle
\big|\frac{f(b)-f(a)}{b-a}-f'(a)\big|\le \frac{1}{k},\,\text{ resp. }
\big|\frac{g(b)-g(a)}{b-a}-g'(a)\big|\le \frac{1}{k}\,.
}
$$ 
Then for any $a\in D(f)\cap D(g)$ and any $b\in M$ with 
$0_K<|a-b|\le\frac{1}{l}$ we have that 
\begin{eqnarray*}
&&{\textstyle
\big|\frac{(x f+y g)(b)-(x f+y g)(a)}{b-a}-(x f'+y g')(a)\big|\le}\\
&&{\textstyle
\le\,|x|\cdot\big|\frac{f(b)-f(a)}{b-a}-f'(a)\big|+
|y|\cdot\big|\frac{g(b)-g(a)}{b-a}-g'(a)\big|\le\frac{|x|+|y|}{k}\,.
}
\end{eqnarray*}
Hence on the set $D(f)\cap D(g)$ the derivative of $xf+yg$ equals $xf'+yg'$ and is uniform.
\eproof

We differentiate products and obtain the countable Leibniz 
formula.

\begin{thm}\label{thm_Leibniz}
We assume that the functions $f,g\in\mathcal{F}(M)$ are bounded, have uniform derivatives and that one of them is {\em UC} and the derivative of the other is bounded. Then
$D(f)\cap D(g)\sus D(fg)\sus M$, the restricted derivative
$$
(fg)'\,|\,D(f)\cap D(g)=f'g+fg'
$$
and is uniform on $D(f)\cap D(g)$. 
\end{thm}
\proof
Let $f$, $g$ and $M$ be as stated and let $f$ be UC and $g'$ be 
bounded, the other case is treated similarly. Let a~$k$ be given.
We take an $l$ such that for any $a\in D(f)$, resp. $a\in D(g)$, and any $b\in M$ with $0_K<|a-b|\le\frac{1}{l}$ it holds that
$${\textstyle
\big|\frac{f(b)-f(a)}{b-a}-f'(a)\big|\le \frac{1}{k},\,\text{ resp. }
\big|\frac{g(b)-g(a)}{b-a}-g'(a)\big|\le \frac{1}{k}\,,
}
$$ 
and that for $a,b\in M$ with 
$|a-b|\le\frac{1}{l}$ always $|f(b)-f(a)|\le \frac{1}{k}$. 
We take a~constant $y$ bounding the sets $f[M]$, $g[M]$ and $g'[D(g)]$. 
Then for any $a$ in $D(f)\cap D(g)$ and any $b\in M$ with 
$0_K<|a-b|\le\frac{1}{l}$ one has that 
\begin{eqnarray*}
&&{\textstyle
\big|\frac{(fg)(b)-(fg)(a)}{b-a}-(f'g+fg')(a)\big|=}\\
&&{\textstyle
=\,\big|\frac{f(b)(g(b)-g(a))+(f(b)-f(a))g(a)}{b-a}-f'(a)g(a)-f(a)g'(a)\big|}\\
&&{\textstyle
\le|f(b)|\cdot\big|\frac{g(b)-g(a)}{b-a}-g'(a)\big|+|f(b)-f(a)|\cdot|g'(a)|+}\\
&&{\textstyle
+\,\big|\frac{f(b)-f(a)}{b-a}-f'(a)\big|\cdot|g(a)|\le\frac{y}{k}+
\frac{y}{k}+\frac{y}{k}=\frac{3y}{k}\,.
}
\end{eqnarray*}
Hence on $D(f)\cap D(g)$ the derivative of $fg$ equals $f'g+fg'$ and is uniform.
\eproof

\noindent
The countable version of the formula for derivatives of 
ratios will be obtained elsewhere, we do not need it here.

We differentiate polynomials.

\begin{prop}\label{prop_deriPoly}
Let $M\sus\Q$ and 
$p\in\mathcal{F}(M)$ 
be a~polynomial. If $p=k_{0_K}$ then the derivative $p'=k_{0_K}\,|\,L(M)$ 
and is uniform. If $p$ has the canonical form 
$p(a)=\sum_{j=0}^n x_ja^j$ then the following hold.
\begin{enumerate}
\item For $n>0$ we have
$p'(a)=\sum_{j=1}^n jx_ja^{j-1}\,|\,L(M)$ and for $n=0$ we have 
$p'=k_{0_K}\,|\,L(M)$.
\item If $M$ is bounded then $p'$ is uniform. 
\item If $\deg p\le 1$ then $p'$ is uniform for any $M$.
\end{enumerate}
\end{prop}
\proof
Let $M$ and $p$ be as stated. The first claim was proven in Proposition~\ref{prop_UderKonIde}.
1 and 2 follow by induction on the degree of $p$ with the help of Propositions~\ref{prop_UderKonIde} and \ref{prop_linDer} and Theorem~\ref{thm_Leibniz}. 3 follows from Propositions~\ref{prop_UderKonIde} and
\ref{prop_linDer}. 
\eproof

The next countable version of VDP (vanishing derivative principle) is 
the most important result of Section~\ref{sec_inteRealF}. 

\begin{thm}\label{thm_GlobEx}
Suppose that $u<x<v$ are in $K$, $I\equiv [u,v]_{\Q}$,   
$f,f'\in\mathcal{F}(I)$, $f'$ is uniform and {\em UC}, and
that $f'(x)\ne0_K$. Then there exist fractions $b',b''\in I$ such that
$$
f(b')<f(x)<f(b'')\;.
$$
Hence the (extended) value $f(x)$ is not a~global extreme of $f$.
\end{thm}
\proof
Let $u$, $x$, $v$, $I$, $x$ and $f$ be as stated. Let $f'(x)<0_K$, the 
case $f'(x)>0_K$ is treated similarly. Note that $f$ is UC by 2 
of Proposition~\ref{prop_UCandDer}. We get the fractions $b'$ and $b''$ 
via the approximation
$$
f(b)=f(a)+\big(f'(a)+\delta_{a,\,f}(b)\big)\cdot(b-a)
$$
where $a,b\in I$ are near $x$. 
We show how to obtain $b'$. 
The obtaining of $b''$ is outlined at the end of the proof.

Since $f'$ is uniform, there is 
an $l$ such that $x+\frac{1}{l}<v$ and 
$${\textstyle
\text{if }a,\,b\in I\wedge a\le b\le a+\frac{1}{l}\,\text{ then }\,|\de_{a,\,f}(b)|\le \frac{|f'(x)|}{3}\,.
}
$$
For every such $a$ and $b$, 
$$
{\textstyle
\text{if also }
f'(a)\le\frac{f'(x)}{2}\wedge b\ge a+\frac{1}{2l}\,\text{ then }\,
\big(f'(a)+\delta_{a,\,f}(b)\big)\cdot(b-a)\le \frac{1}{12l}\cdot f'(x)\,.
}
$$
Recall that $f'(x)<0_K$. We approximate $f(x)$ and $f'(x)$ by, respectively, 
$f(a)$ and $f'(a)$ according to Theorem~\ref{thm_extending} and take an $a>x$ so near to $x$ that 
$$
{\textstyle
a+\frac{1}{l}\le v\wedge f'(a)\le\frac{1}{2}\cdot f'(x)\wedge f(a)+\frac{1}{12l}\cdot f'(x)<f(x)\;.
}
$$
Finally, we take any fraction $b'$ such that $a+\frac{1}{2l}\le b'\le a+\frac{1}{l}$. 
Then the two displayed implications give that indeed
$$
{\textstyle
f(b')=f(a)+\big(f'(a)+\delta_{a,\,f}(b')\big)\cdot(b'-a)\le f(a)+\frac{1}{12l}\cdot f'(x)<f(x)\,.
}
$$

As promised, we explain how to get a~$b''\in[u,v]_{\Q}$ with $f(b'')>f(x)$. We replace the first inequality for $l$ with $u<x-\frac{1}{l}$, replace the next two ranges 
for $b$ with $a-\frac{1}{l}\le b\le a$ and $b\le a-\frac{1}{2l}$, respectively, and reverse the 
final inequality in the second implication to $\ds\ge-
\frac{1}{12l}\cdot f'(x)$. Then we select some $a$ with $a<x$ and
such that $u\le a-\frac{1}{l}$, $f'(a)$ is as before but $f(a)-\frac{1}{12l}\cdot f'(x)>f(x)$. For any $b''$ with $a-\frac{1}{l}\le b''\le a-\frac{1}{2l}$ one has that $f(b'')>f(x)$.
\eproof

\noindent
Since $b'$ and $b''$ may be as close to $x$ as we wish, in the 
theorem we can replace ``global'' with ``local''. In 
Section~\ref{sec_stranFun} we give an example showing that uniformity 
of $f'$ cannot be omitted.

We obtain countable versions of Rolle's and Lagrange's theorems. 

\begin{thm}\label{thm_Rolle}
Suppose that $u<v$ are in $K$, $I\equiv [u,v]_{\Q}$, $f,f'\in\mathcal{F}(I)$, 
$f(u)=f(v)$ and that $f'$ is uniform and {\em UC}. Then there is an $x\in K$ such that
$$
u<x<v\wedge f'(x)=0_K\;.
$$
\end{thm}
\proof
The function $f$ is UC by 2 of Proposition~\ref{prop_UCandDer}. If 
$f$ is constant then $f'(a)=0_K$ for every $a\in I$ by 
Proposition~\ref{prop_UderKonIde} and any $x$ between $u$ and $v$ 
works. 

Suppose that $f(a)>f(u)=f(v)$ for some $a\in I$, the other case with 
$f(a)<f(u)=f(v)$ 
is treated similarly. By Theorem~\ref{thm_attaGlobExtr} there is an $x\in K$ with 
$u\le x\le v$ such that $f(x)$ is a~global maximum of $f$. As $f(x)\ge f(a)$, it follows 
that $u<x<v$. By Theorem~\ref{thm_GlobEx}, $f'(x)=0_K$.
\eproof
\vspace{-3mm}
\begin{thm}\label{thm_Lagrange}
Suppose that $u<v$ are in $K$, $I\equiv [u,v]_{\Q}$, $f,f'\in\mathcal{F}(I)$, and that $f'$ is uniform and {\em UC}. Then there is an $x\in K$ such that
$${\textstyle
u<x<v\wedge f'(x)=\frac{f(v)-f(u)}{v-u}\,.
}
$$
\end{thm}
\proof
The function $f$ is UC by 2 of Proposition~\ref{prop_UCandDer}. 
We set $y\equiv \frac{f(v)-f(u)}{v-u}$ and consider the function $g\in\mathcal{F}(I)$
given by $g(a)\equiv f(a)-ya$.
By Propositions~\ref{prop_UCandSums},
\ref{prop_UderKonIde} and \ref{prop_linDer}, the derivative
$$
g'(a)=f'(a)-y\in\mathcal{F}(I)\,,
$$
is uniform and is UC. Hence $g'(x)=f'(x)-y$ for every $x\in\overline{I}$. Also, $g(v)-g(u)=f(v)-f(u)-y\cdot(v-u)=0_K$, thus $g(u)=g(v)$. We apply 
Theorem~\ref{thm_Rolle} 
to $g$ and get an $x$ such that $u<x<v$ and $g'(x)=f'(x)-y=0_K$. Hence $f'(x)=y=
\frac{f(v)-f(u)}{v-u}$.
\eproof
\vspace{-3mm}
\subsection{The countable exponential function}\label{subsec_expon}

We derive properties of the countable exponential function. 

\begin{defi}\label{def_expFce}
We define the exponential function
$\mathrm{e}^a=\exp\,a
\in\mathcal{F}(\Q)$
as the sum  
$${\textstyle
\exp(a)=\exp\,a\equiv\sum_{n=0}^{\infty}
\frac{a^n}{n!}
=\sum_{n=0}^{\infty}E\big(\frac{a^n}{n!}\big)\,.
}
$$
\end{defi}
Correctness of this definition follows from the absolute 
convergence of the involved series for every $a$, which is easy 
to prove. 

We derive two bounds on $\mathrm{e}^a$ for $a$ near zero, 
and remind Euler's number $\mathrm{e}$.

\begin{prop}\label{prop_expUnuly}
For every $a\in\Q$ with $|a|\le\frac{1}{2}$,  
$$
\big|\mathrm{e}^a-1_K\big|\le 2|a|
\,\text{ and }\,
\big|\mathrm{e}^a-1_K-a\big|\le a^2\,.
$$
\end{prop}
\proof
These bounds follow from Definition~\ref{def_expFce} and Propositions~\ref{prop_linSums},  \ref{prop_ACimplC} and \ref{prop_geomSer}. 
\eproof
\vspace{-3mm}
\begin{defi}\label{def_EulerNum}
The number $\mathrm{e}\equiv\exp\,1$ ($\in K$).
\end{defi}

For the proof of Theorem~\ref{thm_expIden} we need the binomial theorem. We leave 
its combinatorial proof as an exercise for the reader.

\begin{prop}\label{prop_binomT}
Let $x$ and $y$ be real numbers and $n\in\omega$. Then 
$${\textstyle
(x+y)^n=\sum_{i=0}^n\binom{n}{i}\cdot x^i\cdot y^{n-i}\,\text{ where }\,\binom{n}{i}=
\frac{1}{i!}n(n-1)\ds(n-i+1)\,.
}
$$
\end{prop}

The first main property of 
$\mathrm{e}^a$ is the exponential identity.

\begin{thm}\label{thm_expIden}
For every $a,b\in\Q$ the equality $\mathrm{e}^{a+b}=\mathrm{e}^a\cdot\mathrm{e}^b$ holds.
\end{thm}
\proof
Let $a,b\in\Q$. By Theorem~\ref{thm_prodACser}, $\mathrm{e}^a\cdot\mathrm{e}^b=\sum_{n=0}^{\infty}
\frac{a^n}{n!}
\cdot\sum_{n=0}^{\infty}
\frac{b^n}{n!}$ equals
$${\textstyle
\sum_{n=0}^{\infty}\sum_{i=0}^n \frac{a^i}{i!}\cdot\frac{b^{n-i}}{(n-i)!}
=\sum_{n=0}^{\infty}\frac{1}{n!}\sum_{i=0}^n\binom{n}{i}\cdot a^i\cdot b^{n-i}
}
$$
which by Proposition~\ref{prop_binomT}
equals $\sum_{n=0}^{\infty}
\frac{(a+b)^n}{n!}
=\mathrm{e}^{a+b}$.
\eproof
\vspace{-3mm}
\begin{cor}\label{cor_expKladna}
If $a<b$ are in $\Q$ then $0_K<\mathrm{e}^a<\mathrm{e}^b$.
\end{cor}
\proof
Clearly, $\mathrm{e}^0=1_K$ and $\mathrm{e}^a>1_K$ for $a>0$. By Theorem~\ref{thm_expIden}, 
$\mathrm{e}^a\cdot\mathrm{e}^{-a}=\mathrm{e}^0=1_K$ for any~$a$. Thus the exponential function is 
positive. If $a$ and $b$ are as stated then by Theorem~\ref{thm_expIden}, $\mathrm{e}^b-
\mathrm{e}^a=\mathrm{e}^a(\mathrm{e}^{b-a}-1_K)>0_K$.
\eproof

The exponential grows faster than any polynomial.

\begin{cor}\label{cor_limEnaminuxx} For any $m\in\omega$ and any real sequence $(x_n)$ going  to $+\infty$,
$$
\lim_{n\to\infty}\mathrm{e}^{-x_n}x_n^m=0_K\;.
$$
\end{cor}
\proof
Let $m\in\omega$. We may assume that always $x_n\ge0$. 
By Definition~\ref{def_expFce}, $\mathrm{e}^{x_n}\ge x_n^{m+1}\cdot(m+1)!^{-1}$. Thus by Proposition~\ref{prop_limSeqDomi},
$\lim_{n\to\infty}\mathrm{e}^{-x_n}x_n^m=0_K$ since $0_K\le \mathrm{e}^{-x_n}x_n^m\le(m+1)!\cdot x_n^{-1}$.
\eproof

We show that a~slightly generalized exponential is UC on 
bounded sets.

\begin{prop}\label{prop_ExpIsUC}
For every $c\in\Q$ and bounded $M\sus\Q$, $\mathrm{e}^{ca}\,|\,M$ is {\em UC}.
\end{prop}
\proof
Let $c$ and $M$ be as stated and let $m\in\N$  bound $M$. Let a~$k\ge2|c|$ be given and let $a,b\in M$ with 
$|a-b|\le\frac{1}{k}$. By Theorem~\ref{thm_expIden}, Corollary~\ref{cor_expKladna} and 
Proposition~\ref{prop_expUnuly},
$${\textstyle
\big|\mathrm{e}^{ca}-\mathrm{e}^{cb}\big|=
\mathrm{e}^{cb}\cdot\big|\mathrm{e}^{ca-cb}-1_K\big|\le\mathrm{e}^{|cm|}\cdot 2|c|\cdot|a-b|\le 
\frac{2|c|\cdot\mathrm{e}^{|cm|}}{k}\,.
}
$$
Hence $\mathrm{e}^{ca}\,|\,M$ is UC.
\eproof

We differentiate $\mathrm{e}^{ca}$. The second main property of the exponential is the identity $(\mathrm{e}^a)'=\mathrm{e}^a$.

\begin{prop}\label{prop_derExp}
Suppose that $c\in\Q$, $M\sus\Q$ is bounded and $f(a)\equiv 
\mathrm{e}^{ca}\,|\,M$. Then $f'(a)=c\cdot\mathrm{e}^{ca}\,|\,L(M)$, is uniform and {\em UC}.
\end{prop}
\proof
Let $c$, $M$ and $f$ be as stated. Let $d$ bound $M$, $g(a)\equiv c\cdot\mathrm{e}^{ca}$ and let 
a~$k\ge2|c|$ be given. By Proposition~\ref{prop_expUnuly}, 
if $a\in L(M)$ and $b\in M$ with $0_K<|a-b|\le\frac{1}{k}$, then
$${\textstyle
\big|\frac{f(b)-f(a)}{b-a}-g(a)\big|=\big|\mathrm{e}^{ca}\cdot\frac{\mathrm{e}^{cb-ca}-1_K-c(b-a)}{b-a}\big|
\le\mathrm{e}^{|cd|}\cdot|c^2(a-b)|\le
\frac{c^2\cdot\mathrm{e}^{|cd|}}{k}\,.
}
$$
We also used Theorem~\ref{thm_expIden} and 
Corollary~\ref{cor_expKladna}. Hence $f'(a)=g(a)$ on $L(M)$ and is 
uniform. By Propositions~\ref{prop_UCandPr} and \ref{prop_ExpIsUC}, 
$f'$ is UC on $M$.
\eproof

We differentiate the function
$y\cdot a^n\cdot\mathrm{e}^{ca}$.

\begin{prop}\label{prop_polyExpFun}
Suppose that $y\in K$, $n\in\omega$ and $c\in\Q$, that $M\sus\Q$ is bounded and that $f(a)\equiv y\cdot a^n\cdot\mathrm{e}^{ca}\,|\,M$. For $n=0$, the derivative $f'(a)=y\cdot c\cdot\mathrm{e}^{ca}\,|\,L(M)$, is uniform and is {\em UC}. For $n>0$, the derivative
$$
f'(a)=y\cdot(na^{n-1}+ca^n)\cdot\mathrm{e}^{ca}\,|\,L(M)\,,
$$
is uniform and is {\em UC}.
\end{prop}
\proof
This follows from Propositions~\ref{prop_UCandPr}, 
\ref{prop_deriPoly} and \ref{prop_derExp}, Corollary~\ref{cor_XnaNisUC} and Theorem~\ref{thm_Leibniz}.
\eproof

\subsection{Integrals of countable functions}\label{subsec_integr}

Let $u<v$ be in $K$ and $I\equiv [u,v]_{\Q}$. A~{\em partition} of 
$I$ is a~tuple 
$$
\overline{u}=(u_0,\,u_1,\,\ds,\,u_n),\ n\in\N\,,
$$
of $n+1$ real numbers $u_i$ such that $u=u_0<u_1<\ds<u_n=v$. We define 
$\Delta(\overline{u})\equiv \max(\{u_i-u_{i-1}\;|\;i\in[n]\})$ $(>0_K)$. 
A~{\em tagged partition} of $I$ is a~pair $P=(\overline{u},\overline{b})$ 
of a~partition $\overline{u}$ of $I$ and an $n$-tuple $\overline{b}=(b_1,\ds,b_n)$ 
of fractions $b_i$ such that $b_i\in[u_{i-1},u_i]_{\Q}$, $i\in[n]$. We set $\Delta(P)\equiv \Delta(\overline{u})$.
For $f\in\mathcal{F}(I)$ and a~tagged partition $P=(\overline{u},\overline{b})$ 
of $I$ we consider the finite sum
$${\textstyle
R(f,\,P)\equiv \sum_{i=1}^n f(b_i)\cdot (u_i-u_{i-1})\ \ (\in K)\,.
}
$$
If $f$ is UC we allow real tags $\overline{v}=(v_1,\ds,v_n)$ with 
$u_{i-1}\le v_i\le u_i$ because then the values $f(v_i)$ are defined. Real tags appear in Propositions~\ref{prop_shiftInter1} 
and \ref{prop_subsImprInteg}, and in Theorem~\ref{thm_FTAn}. 

For the proof of Theorem~\ref{thm_riemInt} we need a~lemma on two partitions.

\begin{lem}\label{lem_forNextThm}
Suppose that $I=[u,v]_{\Q}$ is as above, $k\in\N$, $(u_0,\ds,u_n)$ and 
$(v_0,\ds,v_m)$ are 
partitions of $I$ such that
$\{u_0,\ds,u_n\}\sus\{v_0,\ds,v_m\}$, 
and that $(x_1,\ds,x_n)$ and $(y_1,\ds,y_m)$ are real tuples 
such that for every $i\in[n]$ and $j\in[m]$,
$${\textstyle
u_{i-1}\le v_{j-1}<v_j\le u_i\Rightarrow|x_i-y_j|\le \frac{1}{k}\,.
}
$$  
Then
$$
{\textstyle
\big|\sum_{i=1}^n x_i\cdot(u_i-u_{i-1})-\sum_{j=1}^m y_j\cdot (v_j-v_{j-1})\big|\le\frac{v-u}{k}\,.
}
$$
\end{lem}
\proof
For $n=1$ this bound follows from the triangle inequality:
\begin{eqnarray*}
&&{\textstyle
\big|x_1\cdot (v-u)-\sum_{j=1}^m y_j\cdot (v_j-v_{j-1})\big|=
\big|\sum_{j=1}^m x_1\cdot (v_j-v_{j-1})\,-}\\
&&{\textstyle
-\,\sum_{j=1}^m y_j\cdot (v_j-v_{j-1})\big|
\le\sum_{j=1}^m|x_1-y_j|\cdot (v_j-v_{j-1})}\\
&&{\textstyle
\le \frac{1}{k}\sum_{j=1}^m (v_j-v_{j-1})=\frac{v-u}{k}\,.}
\end{eqnarray*}
For $n>1$ we use the assumed inclusion and express the stated 
difference as a~sum of differences that we treated for $n=1$; they correspond 
to the partitions of the intervals $[u_{i-1},u_i]_{\Q}$, $i\in[n]$, by 
the fractions $v_j$ lying in them.  By the triangle inequality and the 
case $n=1$ the absolute value of the difference is at most
$\sum_{i=1}^n\frac{u_i-u_{i-1}}{k}=\frac{v-u}{k}$.
\eproof

We introduce integrals of functions in 
$\mathcal{F}([u,v]_{\Q})$.

\begin{defi}\label{def_integrals}
Suppose that $I\equiv [u,v]_{\Q}$ and that $f\in\mathcal{F}(I)$. If for every 
sequence $(P_n)$ of tagged partitions of $I$ with $\lim\Delta(P_n)=0_K$
the limit 
$$
\lim_{n\to\infty} R(f,\,P_n)\ \ (\in K)
$$
exists, we call it the integral of $f$ (from $u$ to $v$) and denote it by
$${\textstyle
(\Q)\int_u^v f\,\text{ or by }\,(\Q)\int_u^v f(a)\,.
}
$$
\end{defi}
The interleaving argument shows 
that the limit, if it exists, is unique and does not depend on the 
sequence $(P_n)$. The integration variable $a$ runs from $u$ to $v$  
through fractions. We set 
$(\Q)\int_u^u f\equiv 0_K$ for any $u$ and any $f$,  and $(\Q)\int_u^v f\equiv -(\Q)\int_v^u f$ 
for $u>v$ if the latter integral exists.

Let $P=(\overline{u},\overline{b})$ and $Q=(\overline{v},\overline{c})$, 
where $\overline{u}=(u_0,\ds,u_l)$ and $\overline{v}=
(v_0\ds,v_m)$, be tagged 
partitions of $I=[u,v]_{\Q}$. Their {\em union $P\cup Q$} is the tagged partition
$$
P\cup Q\equiv (\overline{u}\cup\overline{v},\,\overline{b}\cup\overline{c})
$$ 
of $I$ such that $\overline{u}\cup\overline{v}$ is 
the partition $(w_0,\ds,w_n)$ of $I$, where 
$$
\{w_0,\,\ds,\,w_n\}\equiv 
\{u_0,\,\ds,\,u_l\}\cup\{v_0,\,\ds,\,v_m\}\,,
$$
and $\overline{b}\cup\overline{c}=(d_1,\ds,d_n)$ are arbitrary fractions $d_i
\in[w_{i-1},w_i]_{\Q}$, $i\in[n]$.

We prove the basic existence theorem for the countable integral. 

\begin{thm}\label{thm_riemInt}
If $I\equiv [u,v]_{\Q}$ then
every {\em UC} function $f\in\mathcal{F}(I)$ has $(\Q)\int_u^v f$.
\end{thm}
\proof
Let $I$ and $f$ be as stated, $(P_n)=((\overline{u}_n,\overline{b}_n))$ be 
a~sequence of tagged partitions of $I$ with $\lim \Delta(P_n)=0_K$, and let a~$k$ be 
given. Since $f$ is UC, we take an $l$ such that for any fractions $c$ and $c'$,
$${\textstyle
c,\,c'\in I\wedge|c-c'|\le \frac{1}{l}\Rightarrow|f(c)-f(c')|\le \frac{1}{k}\,.
}
$$
We show that $(R(f,P_n))$ is a~Cauchy sequence. For every 
large $n$, $\Delta(P_n)\le 
\frac{1}{l}$. The previous lemma, the triangle inequality and the 
choice of $l$ imply that 
for every large $m$ and $n$, 
\begin{eqnarray*}
&&|R(f,\,P_m)-R(f,\,P_n)|\le\\
&&\le
|R(f,\,P_m)-R(f,\,P_m\cup P_n)|+|R(f,\,P_m\cup P_n)-R(f,\,P_n)|\le\\
&&\le{\textstyle \frac{v-u}{k}+
\frac{v-u}{k}=\frac{2(v-u)}{k}\,.
}
\end{eqnarray*}
Thus $(R(f,P_n))$ is Cauchy and by Theorem~\ref{thm_complR} it has a~limit.
\eproof

We establish continuity of countable integrals.

\begin{prop}\label{prop_equivDefInteg}
Let $I\equiv [u,v]_{\Q}$ and
$f\in\mathcal{F}(I)$ be {\em UC}. Then for every $k$ there is an $l$ 
such that for every tagged partition $P$ of $I$ (possibly with 
real tags),
$${\textstyle
\Delta(P)\le\frac{1}{l}\Rightarrow\big|R(f,\,P)-
(\Q)\int_u^v f\big|\le\frac{1}{k}\,.
}
$$
\end{prop}
\proof
Let $I$ and $f$ be as stated and let a~$k$ be given. The integral 
exists by Theorem~\ref{thm_riemInt}. We take 
an $l'$ such that 
$\frac{2(v-u)+1_K}{l'}\le\frac{1}{k}$. We use that $f$ is UC and take an $l$ such that
$${\textstyle
u\le x,\,x'\le v\wedge|x-x'|\le\frac{1}{l}\Rightarrow|f(x)-f(x')|\le \frac{1}{l'}\,.
}
$$
We use Definition~\ref{def_integrals} and 
take a~tagged partition $Q$ of $I$ such that 
$\Delta(Q)\le\frac{1}{l}$ and $|R(f,Q)-(\Q)\int_u^v f|\le
\frac{1}{l'}$. Lemma~\ref{lem_forNextThm} and the 
triangle inequality imply that for every tagged partition $P$ of $I$ 
with $\Delta(P)\le \frac{1}{l}$ and possibly real tags,
\begin{eqnarray*}
&&{\textstyle \big|R(f,P)-(\Q)\int_u^v f\big|\le
|R(f,\,P)-R(f,\,P\cup Q)|+|R(f,\,P\cup Q)\,-}\\
&& -\,R(f,\,Q)|\,+\big|R(f,\,Q)-
{\textstyle
(\Q)\int_u^v f\big|\le
\frac{v-u}{l'}+\frac{v-u}{l'}+\frac{1}{l'}
=\frac{2((v-u)+1_K)}{l'}\le
\frac{1}{k}\,.
}
\end{eqnarray*}
\eproof

We establish the linearity and additivity of countable integrals. 

\begin{prop}\label{prop_linearInteg}
Let $I\equiv [u,v]_{\Q}$ and
$f,g\in\mathcal{F}(I)$ be {\em UC} functions. Then 
$${\textstyle
(\Q)\int_u^v(x f+y g)=x\cdot(\Q)\int_u^v f+y\cdot(\Q)\int_u^v g\,.
}
$$
\end{prop}
\proof
The three integrals exist by Proposition~\ref{prop_UCandSums} and 
Theorem~\ref{thm_riemInt}. It is easy to see that for any tagged partition 
$P$ of $I$ it holds that
$$
R(x f+y g,\,P)=x\cdot R(f,\,P)+y\cdot R(g,\,P)\;.
$$
We are done by Definition~\ref{def_integrals} and 
Proposition~\ref{prop_linCombLimSeq}.
\eproof
\vspace{-3mm}
\begin{prop}\label{prop_lineIntegInter}
Let $u<v<w$ be in $K$, $I\equiv [u,w]_{\Q}$ and let $f\in\mathcal{F}(I)$ be {\em UC}. Then 
$${\textstyle
(\Q)\int_u^w f=(\Q)\int_u^v f+(\Q)\int_v^w f\,.
}
$$
\end{prop}
\proof
The three integrals exist by Theorem~\ref{thm_riemInt}.
Any two tagged partitions $P$ and $P'$ of $[u,v]_{\Q}$ and $[v,w]_{\Q}$, respectively, 
merge in the tagged partition $Q$ of $I$ such that
$$
R(f,\,Q)=R(f,\,P)+R(f,\,P')
\,\text{ and }\,\Delta(Q)=
\max(\{\Delta(P),\,\Delta(P')\})\,.
$$
We are done by Definition~\ref{def_integrals} and 
Proposition~\ref{prop_linCombLimSeq}.
\eproof

Shift identities play a~key role in Hilbert's proof of transcendence of $\mathrm{e}$.

\begin{prop}\label{prop_shiftInter1}
Let $u<v$ and $x$ be in $K$, $I\equiv [u+x,v+x]_{\Q}$ and let 
$f\in\mathcal{F}(I)$ be {\em UC}. Then 
$${\textstyle
(\Q)\int_u^v f(a+x)=(\Q)\int_{u+x}^{v+x} f(a)\,.
}
$$
\end{prop}
\proof
Both integrals exist by Theorem~\ref{thm_riemInt}.
Let $P=(\overline{u},\overline{b})$, where 
$\overline{u}=(u_0,\ds,u_n)$ and 
$\overline{b}=(b_1,\ds,b_n)$, be a~tagged partition of $[u,v]_{\Q}$. We shift $P$ by 
$x$ to the tagged partition $P+x\equiv 
(\overline{v},\overline{c})$ of 
$I$, where $\overline{v}=(u_0+x,\ds,u_n+x)$ and 
$\overline{c}=(b_1+x,\ds,b_n+x)$. In general, we get irrational tags $b_i+x$ 
but this is no problem because $f$ is UC. Clearly, 
$\Delta(P)=\Delta(P+x)$ and 
$R(f(a+x),P)=R(f(a),P+x)$.
Thus if  $(P_n)$ is a~sequence of tagged partitions of $[u,v]_{\Q}$ with $\lim \Delta(P_n)=0$ then 
\begin{eqnarray*}
&&{\textstyle
(\Q)\int_u^v f(a+x)=
\lim_{n\to\infty}R(f(a+x),\,P_n)=\lim_{n\to\infty}R(f(a),\,P_n+x)=}\\
&&{\textstyle
=(\Q)\int_{u+x}^{v+x} f(a)\;.}
\end{eqnarray*}
\eproof

We close this subsection with three integral inequalities.

\begin{prop}\label{prop_ineqInteg}
Let $I\equiv [u,v]_{\Q}$ and $f,g\in\mathcal{F}(I)$ be {\em UC} functions such that $f\le g$ on $I$. Then 
$${\textstyle
(\Q)\int_u^v f\le(\Q)\int_u^v g\,.
}
$$
\end{prop}
\proof
The two integrals exist by Theorem~\ref{thm_riemInt}. As $f\le g$ on $I$,
for every tagged partition $P$ of $I$ the inequality $R(f,P)\le R(g,P)$ holds. We are done by Definition~\ref{def_integrals} and
Proposition~\ref{prop_OrderAndLimSeq}.
\eproof
\vspace{-3mm}
\begin{cor}\label{cor_integAbsValIneq}
Let $I\equiv [u,v]_{\Q}$ and $f\in\mathcal{F}(I)$ be {\em UC}. Then 
$${\textstyle
\big|(\Q)\int_u^v f\big|\le(\Q)\int_u^v |f|\,.
}
$$
\end{cor}
\proof
The two integrals exist by Proposition~\ref{prop_UCandAbsVal} 
and Theorem~\ref{thm_riemInt}.
Since $-|f|\le f\le|f|$ on $I$, Propositions~\ref{prop_linearInteg} and \ref{prop_ineqInteg} give the 
inequalities
$$
{\textstyle
-(\Q)\int_u^v|f|\le (\Q)\int_u^v f \le (\Q)\int_u^v|f|}\,.
$$
They are equivalent to the stated inequality.
\eproof
\vspace{-3mm}
\begin{cor}\label{cor_LMbound}
If $I\equiv [u,v]_{\Q}$, $f\in\mathcal{F}(I)$ is {\em UC} and $x\ge0_K$ is a~real number such 
that $|f|\le x$ on $I$, then 
$${\textstyle
\big|(\Q)\int_u^v f\big|\le(v-u)\cdot x\,.
}
$$
\end{cor}
\proof
We apply Proposition~\ref{prop_ineqInteg} on the inequalities $-x\le f$ and $f\le x$ and use the integral  $(\Q)\int_u^v x=(v-u)\cdot x$ which follows from Definition~\ref{def_integrals}.
\eproof

\subsection{Improper integrals of countable functions}\label{subsec_improp}

We introduce integrals of functions in 
$\mathcal{F}([u,+\infty)_{\Q})$.

\begin{defi}\label{def_impropInteg}
Let $I\equiv [u,+\infty)_{\Q}$ and $f\in\mathcal{F}(I)$ be {\em UC} on 
every interval $[u,v]_{\Q}$, $v>u$. If for every real sequence $(x_n)$ 
going to $+\infty$ the limit 
$${\textstyle
\lim_{n\to\infty}(\Q)\int_u^{x_n}f
}
$$
exists, we denote it by 
$${\textstyle
(\Q)\int_u^{+\infty}f\,\text{ or by }\,(\Q)\int_u^{+\infty}f(a)
}
$$
and call it the improper integral of $f$ (from $u$ to $+\infty$).
\end{defi}
The integrals $\int_u^{x_n}f$ exist by Theorem~\ref{thm_riemInt}; for 
$x_n<u$ we define $\int_u^{x_n}f$ arbitrarily. The interleaving argument 
shows that the
limit, if it exists, does not depend on the sequence $(x_n)$. 

Like ordinary integrals, improper integrals enjoy linearity, 
additivity, and shift identities.

\begin{prop}\label{prop_linearInteg2}
Let $I$ and $f,g\in\mathcal{F}(I)$ be as in Definition~\ref{def_impropInteg}. Then the equality
$${\textstyle
(\Q)\int_u^{+\infty}(x f+y g)=x\cdot(\Q)\int_u^{+\infty} f+y\cdot(\Q)\int_u^{+\infty}g
}
$$
holds if the latter two improper integrals exist. 
\end{prop}
\proof
This follows from Propositions~\ref{prop_linCombLimSeq}, \ref{prop_UCandSums} and \ref{prop_linearInteg}.
\eproof

\begin{prop}\label{prop_linInfIntegInter}
Let $I$ and $f\in\mathcal{F}(I)$ be as in Definition~\ref{def_impropInteg} and let $v>u$. Then the equality
$${\textstyle
(\Q)\int_u^{+\infty} f=(\Q)\int_u^v f+(\Q)\int_v^{+\infty}f
}
$$
holds if one of the two improper integrals exists.
\end{prop}
\proof
We denote the former improper integral by $J_1$ and the latter by 
$J_2$. We assume that $J_1$ exists and take any real sequence $(x_n)$ 
going to $+\infty$. We may assume that always $x_n>v$. Using
Propositions~\ref{prop_linCombLimSeq} and \ref{prop_lineIntegInter} we compute

\begin{eqnarray*}
J_2&=&{\textstyle
\lim_{n\to\infty}(\Q)\int_v^{x_n}f=\lim_{n\to\infty}\big(-(\Q)\int_u^v f+(\Q)\int_u^{x_n}f\big)}\\
&=&{\textstyle
-(\Q)\int_u^v f+\lim_{n\to\infty}(\Q)\int_u^{x_n}f
=-(\Q)\int_u^v f+J_1\,, 
}   
\end{eqnarray*}
which is equivalent to the stated equality.
If $J_2$ is assumed to exist, the computation of $J_1$ is similar.
\eproof

\begin{prop}\label{prop_subsImprInteg}
Suppose that $u,x\in K$ and that $I$ and $f\in\mathcal{F}(I)$ are as in Definition~\ref{def_impropInteg}, with $u+x$ in place of $u$. Then the equality
$${\textstyle
(\Q)\int_u^{+\infty}f(a+x)=(\Q)\int_{u+x}^{+\infty}f(a)
}
$$
holds whenever one of the two improper integrals exists.
\end{prop}
\proof
This follows by Proposition~\ref{prop_shiftInter1} because  $(x_n)$ goes to $+\infty$ iff $(x_n+x)$ does. 
\eproof

\subsection[The fundamental theorem of analysis for countable functions]{The FTA for countable functions}

Here is the countable fundamental theorem of analysis (FTA). To prove it, we devised Theorem~\ref{thm_GlobEx} and 
allowed real tags in $R(f,P)$. 

\begin{thm}\label{thm_FTAn}
Let $I\equiv [u,v]_{\Q}$ and $f\in\mathcal{F}(I)$ be such that
$f'\in\mathcal{F}(I)$ and
is uniform and {\em UC}. Then
$${\textstyle
(\Q)\int_u^v f'=f(v)-f(u)\,.
}
$$
\end{thm}
\proof
Let a~$k$ be given. The function $f$ is UC by 2 of 
Proposition~\ref{prop_UCandDer}. The integral exists by 
Theorem~\ref{thm_riemInt}. 
We take a~partition $\overline{u}=(u_0,\ds,u_n)$ of $I$ such that 
$${\textstyle
\big|R(f',\,P)-(\Q)\int_u^v f'\big|\le \frac{1}{k}
}
$$
holds for any tagged partition 
$P=(\overline{u},\overline{v})$ of 
$I$ (i.e., the inequality holds for any $n$-tuple $\overline{v}$ of real numbers 
 in $P$); by Proposition~\ref{prop_equivDefInteg}
it suffices to take any $\overline{u}$ with small enough $\Delta(\overline{u})$. By 
Theorem~\ref{thm_Lagrange} there exist real numbers $x_i$, $i\in[n]$, such that
$$
u_{i-1}<x_i<u_i\,\text{ and }\,f(u_i)-f(u_{i-1})=f'(x_i)\cdot(u_i-u_{i-1})\,.
$$
With $\overline{v}\equiv (x_1,\ds,x_n)$ and $P\equiv 
(\overline{u},\overline{v})$ (tags have to be real, but this is no 
problem as $f'$ is UC) we get that 
$${\textstyle
f(v)-f(u)=\sum_{i=1}^n(f(u_i)-f(u_{i-1}))=\sum_{i=1}^n f'(x_i)(u_i-u_{i-1})=
R(f',\,P)\;.
}
$$
Thus
$${\textstyle
\big|f(v)-f(u)-(\Q)\int_u^v f'\big|=\big|R(f',\,P)-(\Q)\int_u^v f'\big|\le\frac{1}{k}
}
$$
and the stated equality follows by taking $k\to\infty$.
\eproof

The countable FTA takes the following form for improper integrals.  

\begin{cor}\label{cor_FTAimprop}
Suppose that $I\equiv [u,+\infty)_{\Q}$, 
$f,f'\in\mathcal{F}(I)$, that $f'$ is uniform and {\em UC} 
on any rational interval $[u,v]_{\Q}$ with $v>u$, and that
$f$ is such that for any real sequence $(x_n)$ going to $+\infty$ the (unique) limit $L\equiv \lim f(x_n)$ exists. Then
$${\textstyle
(\Q)\int_u^{+\infty}f'=L-f(u)\,.
}
$$
\end{cor}
\proof
$L$ is unique due to the interleaving argument. The rest follows from Definition~\ref{def_impropInteg} and the 
previous theorem.
\eproof

\subsection{Countable Euler's identity}\label{subsec_EulIntAZ}

We arrive at the countable version of Euler's identity. L.~Euler 
discovered the original form in the 18-th century 
(\cite{vara}). Recall that $k!=1\cdot 
2\cdot\ldots\cdot k$ for $k\in\N$, and that $0!=1$.

\begin{thm}\label{thm_Euler}
For every $k\in\omega$,
$${\textstyle
(\Q)\int_{0_K}^{+\infty}a^k\cdot\mathrm{e}^{-a}=k!\,.
}
$$
\end{thm}
\proof
We denote the integrand by $F_k(a)$, the integral by $I_k$ and
proceed by induction on $k$ starting at $k=0$. By Proposition~\ref{prop_polyExpFun}, $(-
\mathrm{e}^{-a})'=F_0$ (on $\Q$). By Proposition~\ref{prop_polyExpFun} and 
Corollaries~\ref{cor_limEnaminuxx} and 
$\ref{cor_FTAimprop}$, $I_0=-0_K+1_K=1_K$. 
Let $k>0$. By Proposition~\ref{prop_polyExpFun}, on $\Q$ it holds that 
$$
F_k'(a)=k\cdot a^{k-1}\cdot\mathrm{e}^{-a}-a^k\cdot\mathrm{e}^{-a}=k\cdot F_{k-1}(a)-F_k(a)\;.
$$ 
Using Proposition~\ref{prop_polyExpFun}, Corollaries~\ref{cor_limEnaminuxx} and $\ref{cor_FTAimprop}$, and induction we get that
$I_k=k\cdot I_{k-1}-(\Q)\int_{0_K}^{+\infty}F_k'=k\cdot(k-1)!-0_K+0_K=k!$.
\eproof

\noindent Proposition~\ref{prop_linearInteg2} and the previous theorem yield the 
following corollary.

\begin{cor}\label{cor_generIden}
For any polynomial $p\in\mathcal{F}(\Q)$ with the canonical form
$p(a)=x_na^n+\ds+x_1a+x_0$,
$${\textstyle
(\Q)\int_{0_K}^{+\infty}p(a)\cdot\mathrm{e}^{-a}=\sum_{i=0}^n x_i\cdot i!\,.
}
$$
\end{cor}

For now we can part with countable univariate real analysis, what we have built suffices 
to prove the transcendence of $\mathrm{e}$ by means of only HMC functions.

\section[Proving transcendence of 
$\mathrm{e}$ by countable functions ]{Transcendence of $\mathrm{e}$ via countable functions}\label{sec_HilbProof}

Recall that $x\in K$ 
is {\em algebraic} if there exist fractions $a_0$, $\ds$, $a_n$, 
$n\in\N$, such 
that $a_n\ne 0$ and 
$\sum_{i=0}^n a_ix^i=0_K$. In other words, $x$ is a~root of a~nonzero 
polynomial with rational coefficients. 
Non-algebraic real numbers are called {\em transcendental}. Transcendence of $\mathrm{e}$ 
was proved first by Ch.~Hermite in \cite{herm} in 1873; \cite{wald} 
discusses his method. 

\begin{thm}[Hermite, 1873]\label{thm_Hermite}
The number $\mathrm{e}=\exp\,1$ is transcendental.    
\end{thm}

Twenty years later D.~Hilbert \cite{hilb} simplified Hermite's proof and it is this proof we adapt here. We assume for the contrary that $\mathrm{e}=\exp\,1$ is algebraic:
$$
a_n\mathrm{e}^n+\ds+a_1\mathrm{e}+a_0=0_K\,,
$$
where $n\in\N$, $a_i\in\Q$ and $a_n\ne0$. We multiply the left-hand side by a~common 
denominator of the $a_i$, take out the maximum 
possible power of $\mathrm{e}$ and get that all $a_i\in\Z$ and 
$a_0\ne0$. We define integral polynomials $p_m\in\mathcal{F}(\Q)$, $m\in\N$, by
$$
p_m(a)\equiv a^m\big((a-1)(a-2)\ds(a-n)\big)^{m+1}
$$
where $n\in\N$ is as in the above displayed equation for $\mathrm{e}$. The next improper integral exists by 
Corollary~\ref{cor_generIden} and splits by Proposition~\ref{prop_linInfIntegInter}:
\begin{eqnarray*}
I(m)&\equiv &{\textstyle
(\Q)\int_{0_K}^{+\infty}p_m(a)\cdot\mathrm{e}^{-a}}\\
&=&{\textstyle(\Q)\int_{0_K}^ip_m(a)\cdot\mathrm{e}^{-a}+
(\Q)\int_i^{+\infty}p_m(a)\cdot\mathrm{e}^{-a}}\\
&\equiv &I_i(m)+J_i(m),\ \ i\in\omega\;.
\end{eqnarray*}
We multiply the above displayed equation for $\mathrm{e}$ by $I(m)$. So for any $m\in\N$,
$$
0_K=\big(a_n\mathrm{e}^n+\ds+a_1\mathrm{e}+a_0\big)I(m)
={\textstyle
\sum_{i=0}^n a_i\mathrm{e}^i\cdot I_i(m)+
\sum_{i=0}^n a_i\mathrm{e}^i\cdot J_i(m)
}\,.
$$
The next corollary is immediate.
\begin{cor}\label{cor_AplusBje0}
For every $m\in\N$, 
$A(m)+B(m)=0_K$, where
\begin{eqnarray*}
A(m)&\equiv&{\textstyle
\sum_{i=0}^n a_i\cdot\mathrm{e}^i\cdot
(\Q)\int_{0_K}^i p_m(a)\cdot\mathrm{e}^{-a}\,\text{ and}
}\\
B(m)&\equiv&{\textstyle
\sum_{i=0}^n a_i\cdot\mathrm{e}^i\cdot
(\Q)\int_i^{+\infty}p_m(a)\cdot\mathrm{e}^{-a}\,.
}
\end{eqnarray*}
\end{cor}
Using the next two propositions, we bring these equalities to contradiction.

\begin{prop}\label{prop_onAm}
For some $y,z>0_K$ and every $m\in\N$, $|A(m)|\le y\cdot z^m$.
\end{prop}
\proof
We set $w\equiv \sum_{i=0}^n |a_i|\cdot\mathrm{e}^i$. Since $0_K<\mathrm{e}^{-a}\le1_K$ for $a\ge0$ 
and $|p_m(a)|\le n^{(n+1)(m+1)}$ for every $a\in[0,n]_{\Q}$, 
Corollary~\ref{cor_LMbound} 
implies that for every $i=0,1,\ds,n$,
$${\textstyle
|I_i(m)|=\big|(\Q)\int_{0_K}^i p_m(a)\cdot\mathrm{e}^{-a}\big|\le n\cdot n^{(n+1)(m+1)}=n^{n+2}\cdot n^{(n+1)m}\;.}
$$
Thus
$${\textstyle
|A(m)|\le w\cdot\sum_{i=0}^n|I_i(m)|\le w\cdot (n+1)n^{n+2}\cdot n^{(n+1)m}\;.}
$$
The stated bound follows, with $y\equiv w\cdot(n+1)\cdot n^{n+2}$ and $z\equiv n^{n+1}$.
\eproof
\vspace{-3mm}
\begin{prop}\label{prop_onBm}
For infinitely many $m\in\N$, $|B(m)|\ge m!$.
\end{prop}
\proof
For $i=0,\ds,n$, we use
Proposition~\ref{prop_linearInteg2} and Theorem~\ref{thm_expIden}, and get that $a_i\cdot\mathrm{e}^i\cdot J_i(m)=a_i\cdot\mathrm{e}^i\cdot(\Q)\int_i^{+\infty} p_m(a)\mathrm{e}^{-a}$ equals to
$$
{\textstyle
a_i\cdot(\Q)\int_i^{+\infty}p_m(a)\cdot\mathrm{e}^i\cdot\mathrm{e}^{-a}=a_i\cdot(\Q)\int_i^{+\infty}p_m(a)\cdot\mathrm{e}^{-(a-i)}\;.
}
$$
By Proposition~\ref{prop_subsImprInteg}
this is
$$
{\textstyle
a_i\cdot(\Q)\int_{0_K}^{+\infty}p_m(a+i)\cdot\mathrm{e}^{-(a-i+i)}=a_i\cdot(\Q)\int_{0_K}^{+\infty}p_m(a+i)\cdot\mathrm{e}^{-a}\;.
}
$$
Corollary~\ref{cor_generIden} and the definition of $p_m(a)$ imply that 
$B(m)=\sum_{i=0}^n a_i\mathrm{e}^i\cdot J_i(m)$ is an integer for every $m$ and is divisible by $m!$. More precisely,
$$
B(m)\equiv a_0(-1)^{n(m+1)}(n!)^{m+1}\cdot m!\ \  (\mathrm{mod}\;(m+1)!)\;.
$$
If $B(m)=0_K$, we divide the congruence 
by $m!$ and deduce that $m+1$ divides $a_0(n!)^{m+1}$. This is impossible if $m+1$ is coprime with 
$a_0\cdot n!$. Since $a_0\ne0$, there exist infinitely many such $m$. Thus $B(m)\ne 0_K$ and $|B(m)|\ge m!$
for infinitely many $m$.
\eproof

We complete the proof of transcendence of $\mathrm{e}$.

\begin{prop}\label{prop_conclusion} 
The assumption that $\mathrm{e}$ is algebraic leads to a~contradiction, 
hence $\mathrm{e}$ is transcendental. 
\end{prop}
\proof
We divide the equality $A(m)+B(m)=0_K$ in Corollary~\ref{cor_AplusBje0} by $m!$ and get by Corollary~\ref{cor_limSeqExpFac}, 
Proposition~\ref{prop_limSeqDomi} and Proposition~\ref{prop_onAm} that
$${\textstyle
\lim_{m\to\infty}\frac{|B(m)|}{m!}=\lim_{m\to\infty}\frac{|A(m)|}{m!}=0_K\;.}
$$
But $\frac{|B(m)|}{m!}\ge1_K$ for infinitely many $m$ by Proposition~\ref{prop_onBm}, and we have a~contradiction.
\eproof

\noindent
This concludes our proof of transcendence of $\mathrm{e}$ using only HMC real 
functions. We followed the proof in \cite{hilb} which uses improper integrals.
The proof in \cite[p.~4]{bake} works with proper integrals $\int_0^t\mathrm{e}^{t-
u}f(u)\,\mathrm{d}u$, where $f(u)$ is a~real polynomial. In \cite{klaz_LWT} we 
discuss three proofs of a~generalization of the transcendence of $\mathrm{e}$, the 
Lindemann--Weierstrass theorem, see 
\cite[pp.~6--8]{bake} for it.

\section{A~counter-intuitive countable real function}\label{sec_stranFun}

In the proof of the next theorem we construct the function mentioned in the Introduction. It shows that in Theorem~\ref{thm_GlobEx} the assumption of uniformity of derivatives cannot be omitted.

A~{\em segment} is any set $S=\{(a,a+w):\;a\in[u,v]_{\Q}\}$ where $w$, $u$ and $v$ are 
real numbers such that $0_K\le u<v\le 1_K$ and, except for the cases $u=0_K$ 
and $v=1_K$, the numbers $u$ and $v$ are irrational. Its {\em left endpoint}, resp. {\em right endpoint}, is the 
point $(u,u+w)$, resp. $(v,v+w)$. Each segment has slope~$1$, so as a~function in $\mathcal{F}([u,v]_{\Q})$ it has
constant derivative $1$. A~point $(u,v)\in K^2$ {\em lies on (the vertical line) $(x=w)$} if $u=w$. 

Let $x>0_K$ and $\overline{u}=(0_K=u_0<u_1<\ds<u_m=1_K)$, $m\in\N$ with $m\ge2$, be a~partition of 
$I\equiv [0_K,1_K]_{\Q}$ such that all $u_i$ with $i\ne0,m$ are irrational and $u_1=\frac{1}{\sqrt{2}}$. An 
{\em $(x,\overline{u})$-saw function $f\in\mathcal{F}(I)$} is the union of segments
$${\textstyle
f=\bigcup_{i=1}^m S_i
}
$$
such that the left endpoint of $S_1$ is $(0_K,0_K)$,  the left, resp. right, endpoint of $S_i$ lies on 
$(x=u_{i-1})$, resp. on $(x=u_i)$, and for $i\ge2$ the right endpoint of $S_{i-1}$ lies by $x$ higher above the left 
endpoint of $S_i$. The vertical segment joining the latter two endpoints is called a~{\em jump}. We say that $f$ is {\em top}
if the right endpoint of $S_1$, the point $(u_1,u_1)=(\frac{1}{\sqrt{2}},\frac{1}{\sqrt{2}})$, lies above all points of $f$.

\begin{thm}\label{thm_nonAvoidUniDer}
Let $I\equiv [0_K,1_K]_{\Q}$. There exists a~{\em UC} function 
$f\in\mathcal{F}(I)$ such that $f$ has the strict global maximum 
$f(\frac{1}{\sqrt{2}})=\frac{1}{\sqrt{2}}$ and $f'=1_K$ on $I$. 
\end{thm}
\proof
We get $f$ as the pointwise limit of a~sequence of functions $f_n\in\mathcal{F}(I)$, $n\in\N$, such that each $f_n$
is a~top $(x_n,\overline{u}_n)$-saw function, where $x_n\equiv\frac{1}{2^{n-
1}\sqrt{2}}$ and the partition $\overline{u}_n$ of 
$I$ has $m=m_n\equiv 2^{n-1}+1$ subintervals. The general form of 
$f_n$ uniquely determines $f_1$ which equals 
$${\textstyle
\{(a,\,a):\;0_K\le a<\frac{1}{\sqrt{2}}\}\cup
\{(a-\frac{1}{\sqrt{2}},\,a-\frac{1}{\sqrt{2}}):\;\frac{1}{\sqrt{2}}< 
a\le1_K\}\equiv S_{1,1}\cup S_{2,1}
}
$$
and is top. No $f_n$ is UC and 
the main task is to ensure that the pointwise limit $f$ is UC. In
$$
{\textstyle
f_n=\bigcup_{i=1}^{m_n}S_{i,\,n}\,,
}
$$
the first leftmost segment $S_{1,n}$ is common for all $f_n$, 
its left and right endpoints are $(0_K,0_K)$ and 
$(\frac{1}{\sqrt{2}},\frac{1}{\sqrt{2}})$, respectively. The 
segments $S_{i,n}$ with 
$i>1$ are not labeled by $i=2,\ds,m_n$ from left to right as 
in the above definition, they get their labels $i$ during the 
construction. In step $n$ of our construction each segment $S_{i,n}$ is divided in 
$S_{i,n}=T_{i,n}\cup U_{i,n}$, where
$T_{i,n}$ is an initial segment and $U_{i,n}$ is a~final segment. In the 
construction the $U$-segments are stable and only get extended. $T_{1,n}\equiv \emptyset$ and $U_{1,n}\equiv 
S_{1,n}$ for every
$n\in\N$. $T_{i,n}\ne\emptyset$ for every $i>1$ and every $n\in\N$.

We begin the construction by ordering the (rational) interval 
$I$ in an injective sequence $(a_m)$, thus 
$I=\{a_m:\;m\in\N\}$. If 
$f_n$ is defined then $I(m,n)$ denotes the unique point 
$(a_m,y)\in f_n$. The construction preserves the following properties. 

\begin{quote}
{\em If $n>n'$ then $U_{i,n}\supset U_{i,n'}$. For every $n\in\N$ and $m\in[n]$, the point $I(m,n)$ lies inside a~segment $U_{i,n}$ of $f_n$, except 
for $a_m=0_K$, resp. $1_K$, when $I(m,n)$ is the left, resp. right, endpoint of $U_{1,n}$, resp. $U_{2,n}$.}
\end{quote}

Let $n=1$. If $I(1,1)\in S_{1,1}=U_{1,1}$, we set 
$T_{2,1}\equiv S_{2,1}$ and $U_{2,1}\equiv \emptyset$. If $I(1,1)\in S_{2,1}$, we take for $U_{2,1}\ne 
S_{2,1}$ any proper final segment in $S_{2,1}$ such that $I(1,1)\in U_{2,1}$ and set 
$T_{2,1}\equiv S_{2,1}\setminus U_{2,1}$ ($\ne\emptyset$).

We describe how $f_{n+1}$ arises from $f_n$. Suppose that $n\ge1$ 
and that $f_n$ is defined, including segments $T_{i,n}$ and $U_{i,n}$, $i\in[m_n]$. Let $I(n+1,n)\in S_{j,n}$ for a~unique 
index $j\in[m_n]$.

{\em The first case is that $I(n+1,n)\in U_{j,n}$}. For every $i\in[m_n]$ we set 
$U_{i,n+1}\equiv U_{i,n}$. 
For $i=2,\ds,m_n$ we take in $S_{i,n}\setminus U_{i,n+1}$ a~nonempty short proper initial segment $T_{i,n}'$ (in a~moment we explain how short it needs to be). For $i=2,\ds, m_n$ we then set
$$
T_{i,\,n+1}\equiv \big(S_{i,n}\setminus U_{i,n+1}\big)\setminus T_{i,n}'\,\text{ and }\,S_{i,\,n+1}\equiv T_{i,\,n+1}\cup U_{i,\,n+1}\;.
$$
For $i=2,\ds, m_n$ we further set
$${\textstyle
T_{i+2^{n-1},\,n+1}=S_{i+2^{n-1},\,n+1}\equiv 
T_{i,\,n}'+
\uparrow\frac{x_n}{2}\,\text{ and }\,
U_{i+2^{n-1},\,n+1}\equiv\emptyset\,.
}
$$
Here $T_{i,n}'+\uparrow\frac{x_n}{2}$ denotes the upward vertical shift 
of the segment $T_{i,n}'$ by $\frac{x_n}{2}$. 
The left endpoint 
of the shifted $T_{i,n}'$ is the midpoint of the jump above the left endpoint of $S_{i,n}$.
$T_{i,n}'$ is so short that the obtained function $f_{n+1}$ is top, after  the 
shifting everything still lies below the line $(y=\frac{1}{\sqrt{2}})$. 

{\em The second case is that $I(n+1,n)\in T_{j,n}$ (then $j\ge2$)}. We proceed as in the first 
case, except that we extend $U_{j,n}$ to a~proper final 
segment $U_{j,n+1}\supset U_{j,n}$ in $S_{j,n}$, long enough so that 
$I(n+1,n)\in U_{j,n+1}$. Then we define $f_{n+1}$ as before.

In both cases it is easy to check that 
$f_{n+1}$ is a~top $(x_{n+1},\overline{u}_{n+1})$-saw function and that the above displayed properties are preserved. We define
$${\textstyle
f\equiv \bigcup_{n=1}^{\infty}\bigcup_{i=1}^{m_n}U_{i,\,n}\;.
}
$$
It follows that
$f\in\mathcal{F}(I)$ and $f'=1_K$ on $I$. From the construction it is also clear that the point 
$(\frac{1}{\sqrt{2}},\frac{1}{\sqrt{2}})$ lies above all pairs (points) in $f$. It remains to prove that $f$ is UC. 

It is clear that there is a~function $n=n(l)\cc\{4,5,\ds\}\to\N$ such that $n(l)$ goes to $+\infty$ and that for every $l\in\{4,5,\ds\}$ it holds for every $i\in[m_{n(l)}]$ that 
$\frac{1}{l}\le|S_{i,n(l)}|_x$, where $|S_{i,n(l)}|_x$ is the length of the projection 
of the segment $S_{i,n(l)}$ on the $x$-axis. 
(We start with $l=4$ because $\frac{1}{4}<1-\frac{1}{\sqrt{2}}<\frac{1}{3}$.)
Now let a~$k$ be given. We take an $l$ large enough such that
$${\textstyle
\frac{1}{l}+2x_{n(l)}=\frac{1}{l}+\frac{2}{2^{n(l)-1}\sqrt{2}}
\le\frac{1}{k}\,.
}
$$
Let $n\equiv n(l)$. For any given $a<b$ in $I$ with $b-a\le \frac{1}{l}$ we estimate  from above the number $|f(a)-f(b)|$. We have 
$a=a_m$ and $b=a_{m'}$ for some $m,m'\in\N$. By the choice of $l$ there is an $i\in[m_n]$ such that $I(m,n)=(a_m,y)\in S_{i,n}$ and 
either $I(m',n)=(a_{m'},y')\in S_{i,n}$ too or $I(m',n)\in S_{i',n}$ where in $f_n$ the segment 
$S_{i',n}$ follows immediately after $S_{i,n}$. In the former case 
$0_K<y'-y\le\frac{1}{l}$ and in the latter case $-x_n<y'-y\le\frac{1}{l}$.
Of course,  $y=f_n(a)$ and $y'=f_n(b)$. By the construction, $0_K\le f(a)-f_n(a)<x_n$ and $0_K\le f(b)-f_n(b)<x_n$. It follows that
$${\textstyle
|f(a)-f(b)|<x_n+\max(\{\frac{1}{l},\,x_n\})<\frac{1}{l}+2x_n\le\frac{1}{k}
}
$$
and $f$ is UC.
\eproof

\section{Concluding remarks}\label{sec_conclRem} 

We plan to extend our investigations of countable real analysis in three 
directions. 

(i) In a~beautiful 
application of complex analysis to number theory,  A.~Baker 
\cite{bake0,bake} effectively solved Thue equations. In 
\cite{klaz_bake} we plan to obtain a~countable version of his solution.

(ii) A~bizarre situation occurs in finite combinatorics  with the standard definition of planar graphs (\cite{dies}). In order to prove according to it that a~graph, for example
$${\textstyle
K_5^-\equiv\big([5],\,\binom{[5]}{2}\setminus\{\{1,\,2\}\}\big)\,,
}
$$
which is always a~hereditarily {\em finite} set, is planar, one 
has to invoke
{\em uncountable} (!) sets, in this example nine noncrossing plane arcs in $\R^2$ 
representing the edges of $K_5^-$. In \cite{klaz_planar} we plan to 
replace the standard definition with a~better one that fixes the 
main drawback. It is not as much the reliance on uncountable sets 
as the fact that, in contrast with many other graph 
classes, the 
standard definition of planar graphs does not (immediately) 
provide any algorithm for deciding planarity of a~given graph.

(iii) The normalized Lebesgue measure $\lambda\cc\Lambda\to[0,1]$ on $[0,1]$,
where $\Lambda\sus\mathcal{P}([0,1])$ is the system of measurable sets, is a~highly uncountable function. Does it have 
a~countable substitute?

\bigskip\noindent
{\em Department of Applied Mathematics\\ 
Faculty of Mathematics and Physics\\
Charles University\\ 
Malostransk\'e n\'am\v est\'\i\ 25\\
118 00 Praha\\
Czechia}\\
{\tt klazar@kam.mff.cuni.cz}
\end{document}